\documentclass[11pt,a4paper,twoside]{article}

\usepackage{graphicx}

\usepackage{amsmath}
\usepackage{latexsym}
\usepackage{amsfonts}
\usepackage{algorithm}
\usepackage{algpseudocode}
\usepackage{comment} 
\usepackage{multirow}

\usepackage{color}
\usepackage{xcolor}
\usepackage[latin1]{inputenc}

\newcommand{\R}{\mathbb{R}}
\newcommand{\N}{\mathbb{N}}

\newcommand{\halmos}{\hfill $\;\;\;\Box$\\}

\newtheorem{lemma}{Lemma}[section]
\newtheorem{theorem}{Theorem}[section]
\newtheorem{example}{Example}[section]

\newtheorem{proposition}{Proposition}[section]
\newtheorem{assump}{Assumption}[section]

\makeatletter
\newenvironment{breakablealgorithm}
  {
   \begin{center}
     \refstepcounter{algorithm}
     \hrule height.8pt depth0pt \kern2pt
     \renewcommand{\caption}[2][\relax]{
       {\raggedright\textbf{\fname@algorithm~\thealgorithm} ##2\par}%
       \ifx\relax##1\relax 
         \addcontentsline{loa}{algorithm}{\protect\numberline{\thealgorithm}##2}%
       \else 
         \addcontentsline{loa}{algorithm}{\protect\numberline{\thealgorithm}##1}%
       \fi
       \kern2pt\hrule\kern2pt
     }
  }{
     \kern2pt\hrule\relax
   \end{center}
  }
\makeatother

\usepackage{comment}
\usepackage{subcaption}

\makeatletter
\newcommand{\dummylabel}[2]{\def\@currentlabel{#2}\label{#1}}
\makeatother

\begin{document}

\title{A Jacobi-type Newton method for Nash equilibrium problems with descent guarantees}

\author{L. F. Bueno\thanks{Institute of Science and Technology, Federal
    University of S\~ao Paulo, S\~ao Jos\'e dos Campos, SP,
    Brazil. {\tt (lfelipebueno@gmail.com)}} \and G. Haeser\thanks{Department of Applied Mathematics, University of S\~ao Paulo, S\~ao
    Paulo, SP, Brazil. {\tt (ghaeser@ime.usp.br)}} \and O. Kolossoski\thanks{Instituto Latino-Americano de Ci\^encias da Vida e Natureza, Federal University for Latin American Integration, Foz do Igua\c cu, PR, Brazil. {\tt (OliverKolossoski@gmail.com)}}}

\maketitle

\begin{abstract}
A common strategy for solving an unconstrained two-player Nash equilibrium problem with continuous variables is applying Newton's method to the system obtained by the corresponding first-order necessary optimality conditions. However, when taking into account the game dynamics, it is not clear what is the goal of each player when considering they are taking their current decision following Newton's iterates. In this paper we provide an interpretation for Newton's iterate  as follows: instead of minimizing the quadratic approximation of the objective functions parameterized by the other player current decision (the Jacobi-type strategy), we show that the Newton iterate follows this approach but with the objective function parameterized by a prediction of the other player action. This interpretation allows us to present a new Newtonian algorithm where a backtracking procedure is introduced in order to guarantee that the computed Newtonian directions, for each player, are descent directions for the corresponding parameterized functions. Thus, besides favoring global convergence, our algorithm also favors true minimizers instead of maximizers or saddle points, unlike the standard Newton method, which does not consider the minimization structure of the problem in the non-convex case. Thus, our method is more robust compared to other Jacobi-type strategies or the pure Newtonian approach, which is corroborated by our numerical experiments. We also present a proof of the well-definiteness of the algorithm under some standard assumptions, together with a preliminary analysis of its convergence properties taking into account the game dynamics.
\end{abstract}

\newpage

{bf MSC} 49K99, 65K05, 90C30, 91A99

\section{Introduction}

The Nash Equilibrium Problem (NEP), introduced in \cite{nash}, models an economic game where each player aims at minimizing their own objective function, but the decision taken by each player influences the payoff of the other players.  This is a fundamental problem in economics and social behavior theory that has been studied
extensively. See, for instance, \cite{debreu}, \cite{myerson} and \cite{aubin}. For surveys on the subject, see \cite{facchinei} and \cite{fischer}. 
In this paper we consider the unconstrained  two-player NEP given by

\begin{minipage}{0.45\textwidth}
\begin{equation} \label{PJ1}
\begin{array}{cl}
\text{Minimize }& f_1(x_1,x_2),\\[-1pt]
\resizebox{!}{0.25cm}{$x_1 \in  \mathbb{R}^{n_1}$}\\[4pt]
\end{array}
\end{equation}
\end{minipage}
\begin{minipage}{0.45\textwidth}
\begin{equation*}
\begin{array}{cl}
\text{Minimize } & f_2(x_1,x_2),\\[-1pt]
\resizebox{!}{0.25cm}{$x_2 \in  \mathbb{R}^{n_2}$}\\[4pt]
\end{array}
\end{equation*} \end{minipage}\\
where $f_1,f_2\colon\R^{n_1}\times\R^{n_2}\to\R$ are (possibly non-convex) twice continuously differentiable functions that describe each player's objective to be minimized, parameterized by the other player's decision. A (local) solution to the NEP \eqref{PJ1} is a point $(x_1^*,x_2^*)\in\R^{n_1}\times\R^{n_2}$ such that $x_1^*$ is a (local) minimizer for the first problem in \eqref{PJ1} with $x_2=x_2^*$, and $x_2^*$ is a (local) minimizer for the second problem in \eqref{PJ1} with $x_1=x_1^*$.

In order to solve the NEP \eqref{PJ1}, some studies deal with a reformulation of the problem, such as a variational inequality problem or as a complementarity problem; see for instance \cite{majig} and \cite{facchinei-vi}. One problem with this approach is that several solutions of the reformulated problem may not be solutions of the associated equilibrium problem. Another possibility is to reformulate the NEP as an optimization problem via the Nikaido-Isoda function (see \cite{nikaidoisoda}), which transforms the NEP into a minimax problem. A penalty update scheme was also proposed in the works of \cite{fukushima} and \cite{facchinei2009}. Those indirect approaches often require solving a nontrivial optimization problem at each step of the algorithm, an exception being \cite{gradiente} which approximates the Nikaido-Isoda function by replacing the minimization problem by a Cauchy step. In general Nikaido-Isoda based formulations deal with convex objective functions since otherwise the subproblem that arises could be hard to deal with. There are some exceptions, such as \cite{wachter}, which requires pseudoconvexity of a regularized Nikaido-Isoda operator, or \cite{simone}, which presents a branch and bound algorithm on the case of discrete constrained non-convex games. Recent advances deal with the subproblems via similar branch and bound schemes (such as \cite{kirst}) or genetic algorithms (see \cite{indianos}), thus being able to treat a broader class of problems, such as general non-convex or non-differentiable games; even so, the subproblems that arise can still be hard to solve. In order to avoid this kind of computations, the method proposed in this paper aims at solving the system \eqref{PJ1} without rewriting it as an indirect problem, and as we shall see, this leads to solving a linear system of equations at each iteration, while still being able to deal with the non-convex case.

Another recent trend in the literature is considering reinforced learning, see for instance \cite{ozdaglar}, \cite{gamelearning} or \cite{cesarbianchi}. These methods came up from the analysis of equilibrium problems with discrete strategy sets, where the actions of each player are taken considering an aggregate of the other player's previous strategies instead of only their last action. However, expanding on this idea considering non-discrete strategy sets often leads to control problems such as \cite{bervoets}, which may be hard to solve.

One may also consider partitioning methods, namely Jacobi and Gauss-Seidel-type splitting methods (see \cite{Brown}, \cite{facchinei} and \cite{yuan}) that deal with problem \eqref{PJ1} by iterating the idea of fixing the variables of one player and taking a step towards minimizing the problem of the other player. In particular, a standard Jacobi strategy would fix preliminary decisions for each player, and then iteratively update them simultaneously by solving a model for the corresponding optimization problem parameterized by the fixed decision of the other player. Alternatively, so-called best response methods fix the decision of each player as the best possible response and  are commonly extended with relaxation strategies, which expand the range of problems on which the method works and improve its convergence rates, see \cite{caruso-velho}, \cite{caruso}, and \cite{Basar}. The advantage of these approaches is that they are able to treat problem \eqref{PJ1} in its original form, without losing information due to reformulations. However, these methods also require solving a general optimization problem at each step for finding the best response for one player, which may also be computationally expensive, in addition to depending on convexity of the objective functions as conditions for guaranteeing convergence.

In this paper we introduce a method based on a Jacobi-type approach within a Newtonian framework. The main novelty of our approach is that instead of taking a decision based on the current decision of the opponent, our Newton model is built around a prediction of the opponent step. This may reflect the reality of some equilibrium games in practice, in situations where each player, knowing a previous decision of the other player, try to make the best decision based on what is the expected response of the other player. The algorithm proposed in this paper is posed in such a way that each player aims at decreasing their objective function, with respect to the predicted decision of the other player, using quadratic approximations. This avoids the heavy computations of best response methods, while possibly still reflecting the dynamics of equilibrium games in practice.

Additionally, we note that our method considers the true minimization structure of the system \eqref{PJ1} in a Newtonian framework, instead of considering only the system of equations given by the necessary optimality conditions associated with each problem. The idea is that one should favor minimizers in place of maximizers or saddle points when considering the solutions of this system of equations, which is usually not considered in other Newtonian algorithms for solving non-convex instances of \eqref{PJ1}.

In a related work,  a Jacobi-type approach using  quadratic approximations of  players' objective functions was considered in \cite{yuan}. The proposed algorithm consists in a trust region method  which does not directly use the Hessian information on its iterative step; however, it relies on some strong assumptions on the true Hessians that we do not rely on in our work. Additionally, for each player, the method in \cite{yuan} is based on the current response of the other players, rather than on a predicted one.

The paper is organized as follows: In Section \ref{sec2} we introduce our Jacobi-type algorithm in a Newtonian framework.  In Section \ref{boadef} we establish the well-definiteness of our algorithm. In Section \ref{sec3} we state some theoretical properties of the algorithm, establishing convergence to stationarity under certain conditions. 
In section \ref{sec4} we illustrate the behavior and effectiveness of the proposed method with some numerical experiments. Finally, some concluding remarks and future prospective works are given in Section \ref{sec5}.

{\bf Notation:} Given a twice continuously differentiable function $f\colon\R^{n_1+n_2}\to\R,$\linebreak $(x_1,x_2)\mapsto f(x_1,x_2)$, for $i=1,2$ we denote by $\nabla_{x_i}f(\overline{x})$ its partial gradient with respect to variable $x_i\in\R^{n_i}$ evaluated at $\overline{x}\in\R^{n_1+n_2}$. The partial Hessian with respect to variables $x_i$ and $x_j$ evaluated at $\overline{x}$ is denoted by $\nabla_{x_ix_j}^2f(\overline{x})$, for $i,j=1,2$.  Additionally, for each player $i=1,2$, we denote the index of the other player by $\neg i$. We use $\|\cdot\|$ to denote the $2$-norm of vectors and matrices, while $u^Tv$ denotes the canonical inner product of vectors $u$ and $v$ of the same dimensions. We denote by $\mathbb{N}$ the set of positive integers. Given $p \in \N$,  $x \in \R^p$ and $\delta>0$, we denote by $B(x,\delta)=\{y \in \R^p: \|y-x\| \leq \delta\}$, the closed ball of $\R^p$ with center $x$ and radius $\delta$.

\section{A Descent Newton Algorithm for the Two-Player Nash Equilibrium Problem}\label{sec2} 

A well-studied class of algorithms for solving the equilibrium problem \eqref{PJ1} is known as best response methods, whose foundation dates back to the work of \cite{Brown}. In this type of algorithm, it is assumed that, from an initial configuration, a solution is approximated by a dynamic process in which each player makes a decision reacting to the decision of the other player. Thus, it is generally assumed that the functions $b_1(x_2)\colon \mathbb{R}^{n_2} \to \mathbb{R}^{n_1}$ and $b_2(x_1)\colon \mathbb{R}^{n_1} \to \mathbb{R}^{n_2}$ such that
\begin{equation}\label{melhorresposta}b_i(x_{\neg i}) = argmin_{x_i \in \mathbb{R}^{n_i}} f_i(x_i,x_{\neg i})
\end{equation}
are available. Among the algorithms in this class, the most classical ones are of the Jacobi-type, also known as the Cournot method (introduced in \cite{cournot}), where the sequence ${(x_1^k,x_2^k)}$ is defined as $x_i^{k+1}= b_i(x_{\neg i}^k)$, and of the Gauss-Seidel-type, where $x_1^{k+1}= b_1(x_2^k)$ and $x_2^{k+1}= b_2(x_1^{k+1})$;  Classic studies on frameworks related to Jacobi and Gauss-Seidel formulations can be found in \cite{facchineidecomposition}.
More recently, Morgan et. al. developed further the theory on best response methods considering different hypotheses for the objective functions and on different contexts, see \cite{caruso}, \cite{ref2} and references therein. Leslie et al. \cite{leslie} generalized best-response dynamics to stochastic games. In \cite{baudinlaraki} there is an extention of such work for ficticious plays. In \cite{leishanbhag} it is studied bounds on the convergence rate of stochastic equilibrium problems, and in \cite{bayer} sufficient conditions for convergence for best-response dynamics on the context of network games are given. In \cite{passacantando} is devised a best response method for a particular class of functions which yields a sequence of Newtonian linear systems. 


On the other hand, recall that a solution $(x_1,x_2)$ to the equilibrium problem \eqref{PJ1} satisfies the following first-order necessary optimality condition:
 \begin{equation} 
 \label{sistKKT}
\left(\begin{array}{c}
\nabla_{x_1} f_1(x_1,x_2)\\[3pt]
\nabla_{x_2} f_2(x_1,x_2)
\end{array} \right)=
\left(\begin{array}{c}
0\\[3pt]
0
\end{array} \right).
\end{equation}

Most algorithms for solving the NEP \eqref{PJ1} rely on approaching directly the nonlinear system of equations \eqref{sistKKT}, disregarding the minimization structure of problem \eqref{PJ1}; see, for instance, \cite{fischer} and the references therein. This poses no issues when $f_1$ and $f_2$ are convex with respect to their decision variables, however, in the presence of non-convexities, the algorithm may find local maximizers or saddle points as often as local minimizers. Our goal in this paper is to combine the ideas of Jacobi and Newton methods in order to propose an Newtonian algorithm that takes into account the minimization structure of the problem in order to find local minimizers more often, thus better reflecting the behavior of two-player games in practice.

Given a current iterate $(x_1^k,x_2^k)$, the classical Newtonian step $(d_1^k,d_2^k)$ for solving \eqref{sistKKT} is given by solving the following linear system:
\begin{equation}\label{std} \left(
\begin{array}{cc}
\nabla^2_{x_1x_1} f_1(x_1^k,x_2^k) & \nabla^2_{x_1x_2} f_1(x_1^k,x_2^k) \\[3pt]
 \nabla^2_{x_2x_1} f_2(x_1^k,x_2^k)&\nabla^2_{x_2x_2} f_2(x_1^k,x_2^k)
\end{array} \right)
\left(\begin{array}{c}
d_1^k\\
d_2^k
\end{array} \right)=-
\left(\begin{array}{c}
\nabla_{x_1} f_1(x_1^k,x_2^k)\\[3pt]
\nabla_{x_2} f_2(x_1^k,x_2^k)
\end{array} \right).
\end{equation}

It is well known that, for optimization problems, applying Newton's method to the corresponding first-order necessary optimality conditions can be better interpreted as minimizing a quadratic approximation of the objective function. With this interpretation, we can extend Newton's method to non-convex problems, 
and we can also  provide elements for globalization schemes. Two well-established strategies for doing this are the use of trust regions and line-search techniques, ensuring progress near the current point. The first one guarantees the existence of a minimizer for the quadratic approximation of the objective function by establishing a compact domain for the subproblem. In line-search approaches, one should use some positive definite approximation for the Hessian of the objective function to obtain a coercive model and descent directions. In this case, one must carefully choose a convenient step-size along the Newton direction.

In order to exploit the minimization structure of the NEP \eqref{PJ1}, we shall consider an idea which resembles a standard line-search  globalization  described above. Given an approximate solution $(x_1^k,x_2^k)$ at some iteration $k$, the most natural step to be made for each player  would be associated with minimizing, with respect to variable $d_i$, a simplified model for the function $f_i(x_i^k+d_i,x_{\neg i}^k)$: 
 \begin{equation}\label{fjac1}
\begin{array}{cl}
\text{Minimize} & \phi_i(x_i^k+d_i,x_{\neg i}^k),\\[-5pt]
\resizebox{!}{0.18cm}{$d_i$}\\[2pt]
\end{array}
\end{equation}
where $\phi_i$ is an approximation of $f_i$.
 
This approach is related to Jacobi's method and its theoretical investigation has been conducted in \cite{yuan} under a trust-region globalization. The main concern that we want to address is to base the minimization model for each player on a prediction of the other player's action, rather than on its actual decision. Thus, 
 the main contribution of this work is considering the case when the direction $d_1^k$, the decision of the first player, is obtained by a Newton step with respect to a predicted behavior of the second player. Namely, given $x_2^k+d_2^k$, a predicted decision for player two, the problem we aim to solve for player one is:
\begin{equation}\label{ouraprox}
\begin{array}{cl}
\text{Minimize} & \phi_1(x_1^k+d_1,x_2^k+d_2^k).\\[-5pt]
\resizebox{!}{0.18cm}{$d_1$}\\[2pt]
\end{array}
\end{equation}

Inspired by a Newtonian appeal, let us consider the  following convex quadratic approximation:
\begin{equation}\label{auxeq00}
f_1(x_1^k+d_1,x_2^k+d_2^k)\approx f_1(x_1^k,x_2^k+d_2^k)+\nabla_{x_1} f_1(x_1^k,x_2^k+d_2^k)^T d_1+\frac{1}{2}d_1^TH_1^kd_1,\end{equation}
where $H_1^k$
is a positive definite approximation for the Hessian $\nabla^2_{x_1x_1}f_1(x_1^k,x_2^k+d_2^k)$. Since $d_2^k$ needs to be found together with $d_1^k$, suppose that the same approach is done simultaneously by the other player in order to make their decision $d_2^k$. Then, given $(x_1^k,x_2^k)$, we would still have a Nash equilibrium problem in the variables $(d_1,d_2)$. The resulting NEP is simpler to solve, since the objective functions are parameterized convex quadratics in the corresponding decision variables. But, since the term $\nabla_{x_1} f_1(x_1^k,x_2^k+d_2^k)^Td_1$ might still combine the variables in a non-linear way, the approximated problem might still be hard to solve. Thinking in eliminating this inconvenience with a typical idea from Newton's method, we use linear approximations of the gradients in the following way:
\begin{equation}
\nabla_{x_1} f_1(x_1^k,x_2^k+d_2^k) \approx \nabla_{x_1} f_1(x_1^k,x_2^k)+ \nabla^2_{x_1x_2} f_1(x_1^k,x_2^k)d_2^k.
\label{aproxgrad}
\end{equation}

In this way, our approximated problem \eqref{ouraprox} is now clearly defined by combining approximations \eqref{auxeq00} and \eqref{aproxgrad}, which has a simple structure and its solution is readily given by solving the linear system:
$$H_1^kd_1 + \nabla^2_{x_1x_2} f_1(x_1^k,x_2^k)d_2^k+\nabla_{x_1} f_1(x_1^k,x_2^k) =0.$$

Repeating the approach for the second player, with a  positive definite $H_2^k \approx \nabla^2_{x_2x_2} f_2(x_1^k+d_1^k,x_2^k)$, our iteration is based on the solution of the following linear system:
  \begin{equation} \left(
\begin{array}{cc}
H_1^k & \nabla^2_{x_1x_2} f_1(x_1^k,x_2^k) \\[3pt]
 \nabla^2_{x_2x_1} f_2(x_1^k,x_2^k)&H_2^k  
\end{array} \right)
\left(\begin{array}{c}
d_1^{}\\
d_2^{}
\end{array} \right)=-
\left(\begin{array}{c}
\nabla_{x_1}f_1(x_1^k,x_2^k)\\[3pt]
\nabla_{x_2}f_2(x_1^k,x_2^k)
\end{array} \right).
\label{NewtonsistKKT}
\end{equation}

Given the equivalence with the standard Newtonian system \eqref{std}, it is natural to choose $H_1^k \approx \nabla^2_{x_1x_1} f_1(x_1^k,x_2^k)$ and $H_2^k \approx \nabla^2_{x_2x_2} f_2(x_1^k,x_2^k)$. Thus, we have associated the standard Newtonian system \eqref{std} with the original NEP, where every player makes its decision considering the predicted decision of the other player. It now seems natural to consider descent conditions associated with the predicted functions $f_1(x_1,x_2^k+d_2^k)$ and $f_2(x_1^k+d_1^k,x_2)$. However, we cannot guarantee that the resulting directions from \eqref{NewtonsistKKT} are descent directions for these functions, as is shown in the following example.
\begin{example}
\label{nodescent}
Consider the equilibrium problem where $f_1(x_1,x_2) := \frac{1}{2} x_1^2+(x_2^2+2x_2+1)x_1$ and
$f_2(x_1,x_2):= \frac{1}{2} (x_2+2)^2$. In this case, we have 
\begin{itemize}
	\item $\nabla_{x_1} f_1(x_1,x_2)=x_1+x_2^2+2x_2+1$, $\nabla^2_{x_1x_1} f_1(x_1,x_2)=1$, 
	$\nabla^2_{x_1x_2} f_1(x_1,x_2)=2 x_2+2$;
	\item $\nabla_{x_2} f_2(x_1,x_2)=x_2+2$, $\nabla^2_{x_2x_2} f_2(x_2,x_2)=1$, and $\nabla^2_{x_2x_1} f_2(x_2,x_2)=0$. 
\end{itemize}

Taking $(x_1^k,x_2^k):=(0,0)$, the Newton direction $(d_1^k,d_2^k)$ obtained from \eqref{NewtonsistKKT} is the solution of:
\begin{equation*} \left(
\begin{array}{cc}
1 & 2 \\[3pt]
0&1
\end{array} \right)
\left(\begin{array}{c}
d_1^{}\\
d_2^{}
\end{array} \right)=-
\left(\begin{array}{c}
1\\[3pt]
2
\end{array} \right).
\end{equation*}
So $d_1^k=3$ and $d_2^k=-2$. Note that $\nabla_{x_1} f_1(x_1^k,x_2^k)=\nabla_{x_1} f_1(x_1^k,x_2^k+d_2^k)=1$ and therefore $d_1^k$ is not a descent direction for  $f_1(x_1,x_2^k+d_2^k)$ (nor for $f_1(x_1,x_2^k)$) at $x_1^k$.	
\end{example}
This means that we do not expect the descent condition $f_1(x_1^k+td_1^k,x_2^k+d_2^k)<f_1(x_1^k,x_2^k+d_2^k)$ to hold even for arbitrarily small $t$. Therefore, an algorithm that requires an Armijo-like line-search could be not well-defined. The reason for not obtaining a descent direction when minimizing the convex quadratic approximation of the function comes from the fact that we did not use its true gradient. That is, the approximation of $\nabla_{x_1} f_1(x_1^k,x_2^k+d_2^k)$ in \eqref{aproxgrad} might not be accurate enough when $\|d_2^k\|$ is large.
To improve this approximation, we should consider also taking small steps along the direction $d_2^k$. In this way, the line-search should be done simultaneously for $f_1$ and $f_2$, considering a simultaneous backtracking 
in $d_1^k$ and $d_2^k$. Thus, replacing $d_2^{k}$ by $td_2^{k}$ in the previous derivation of \eqref{NewtonsistKKT}, we arrive at the following system of equations for computing $(d_1^{},d_2^{})$, given a tentative step length $t$:
\begin{equation*} \left(
\begin{array}{cc}
H_1^k & t\nabla^2_{x_1x_2} f_1(x_1^k,x_2^k) \\[3pt]
 t\nabla^2_{x_2x_1} f_2(x_1^k,x_2^k)&H_2^k  
\end{array} \right)
\left(\begin{array}{c}
d_1^{}\\
d_2^{}
\end{array} \right)=-
\left(\begin{array}{c}
\nabla_{x_1}f_1(x_1^k,x_2^k)\\[3pt]
\nabla_{x_2}f_2(x_1^k,x_2^k)
\end{array} \right).
\end{equation*}

Note that for $t>0$ sufficiently small, we expect the solution $(d_1^{},d_2^{})$ to be such that $d_1^{}$ is a descent direction for the first problem in \eqref{PJ1} while $d_2^{}$ is a descent direction for the second problem in \eqref{PJ1} considering the predicted decision for the other player. Hence, our method favors local minimizers for solving problem \eqref{PJ1}. Notice, however, that the solution $(d_1^{},d_2^{})$ depends on the choice of the step length $t$, thus, once a step length is rejected by the descent condition, a new direction $(d_1^{},d_2^{})$ must be recomputed, similarly to a trust region approach. Finally, once the descent condition is met for some $t$, the new iterate is defined as $(x_1^{k+1},x_2^{k+1}):=(x_1^k,x_2^k)+t(d_1^{},d_2^{})$.

Our algorithm selects the step-sizes based on a standard Armijo-type condition, namely, checking whether the decrease along a component (with the other component fixed at the predicted point) is proportional to what is predicted by the first order approximation of the function. In mathematical terms, for some $\alpha\in(0,1)$,  we check for each $i\in\{1,2\}$ if
\begin{equation*}
f_i(x_i^k+td_i^{},x_{\neg i}^k+td_{\neg i}^{}) \leq f_i(x_i^k,x_{\neg i}^k+td_{\neg i}^{})+\alpha t \nabla_{x_i}f_i(x_i^k,x_{\neg i}^k+td_{\neg i}^{})^Td_i^{},
\end{equation*}
decreasing $t$ and recomputing the directions whenever one of these inequalities do not hold.

Before formally presenting the algorithm, let us discuss some other conditions needed in the line-search procedure. Differently from  optimization problems, since the directions are not necessarily descent directions when $t$ is large, we must ensure a negative slope over the search direction when computing the step-size $t$. We also avoid directions that are near orthogonal to the gradient of the predicted function. Additionally, in order to avoid stagnation when far from stationary points, we also demand that the direction is not too small compared with the corresponding gradient of the predicted function. Therefore, for $\theta \in (0,1)$ and $\gamma>0$, our line-search checks also if 
\begin{align}
 \nabla_{x_1}f_1(x_1^k,x_2^k+td_2^{})^Td_1^{} &\leq -\theta\|\nabla_{x_1}f_1(x_1^k,x_2^k+td_2^{})\|\cdot\|d_1^{}\|, \nonumber \\
 \|d_1^{}\| &\geq \gamma\|\nabla_{x_1}f_1(x_1^k,x_2^k+td_2^{})\| \label{danterior},
\end{align}
with an analogous requirement for the second player. Notice that, in optimization, similar inequalities are standard requirements on the direction for guaranteeing global convergence of a descent method using an Armijo line-search.

Guided by all these ideas, we formalize below the algorithm for computing a sequence $\{(x_1^{k},x_2^{k})\}_{k\in\mathbb{N}}$ for solving the NEP \eqref{PJ1}, where we note that the algorithm includes a safeguarding strategy that replaces the mixed hessians with zeroes whenever stationarity for the corresponding player has been reached and $t$ is small enough, which will be explained later.\\

\begin{breakablealgorithm}
	\caption{Jacobi-type descent Newton algorithm}\label{Algorithm 1}
	\begin{description}
		\item[Step 0.] Given $(x_1^0,x_2^0) \in \R^{n_1}\times\R^{n_2}$, let $\alpha \in (0,1)$, $\tau \in (0,1]$, $\theta>0, \gamma > 0$ and $k:=0$.
		\item[Step 1.] Compute $g_1^k:=\nabla_{x_1}f_1(x_1^k,x_2^k)$ and $g_2^k:=\nabla_{x_2}f_2(x_1^k,x_2^k)$, choose symmetric positive definite matrices $H_1^k$ and $H_2^k$, and set $t:=1$. 
			\item[Step 2.] For each $i \in \{1,2\}$, set
\begin{equation}\label{separar}
					 M_i := \left\{ \begin{array}{cc}\nabla^2_{x_ix_{\neg i}} f_i(x^k), & \text{ if } \|g_i^k\| > 0 \text{ or } t>\tau,
						\\ 0, & \text{ otherwise.}
                        \end{array}\right.
                        \end{equation}
			
		\item[Step 3.] Check if the matrix 
			$$\left[\begin{array}{cc}H_1^k & tM_1 \\[3pt]
				tM_2&H_2^k\end{array}\right]$$
			is non-singular. If not, set $t:=\frac{t}{2}$ and repeat Step 3.
			\item[Step 4.] 
			Find $(d_1^{},d_2^{})$ by solving the linear system
			\begin{equation} \left(
				\begin{array}{cc}
					H_1^k & tM_1 \\[3pt]
					tM_2 &H_2^k  
				\end{array} \right)
				\left(\begin{array}{c}
					d_1^{}\\
					d_2^{}
				\end{array} \right)=-
				\left(\begin{array}{c}
					g_1^k\\[3pt]
					g_2^k
				\end{array} \right).
				\label{NewtonsistKKT2}
			\end{equation}
			\item[Step 5.]  Check the inequalities
				\begin{align}
					 f_1(x^k+td) \leq f_1(x_1^k,x_2^k+td_2^{})+\alpha t \nabla_{x_1}f_1(x_1^k,x_2^k+td_2^{})^Td_1^{}, \label{armijo1}\\ 
					 \nabla_{x_1}f_1(x_1^k,x_2^k+td_2^{})^Td_1^{} \leq -\theta\|\nabla_{x_1}f_1(x_1^k,x_2^k+td_2^{})\|\cdot\|d_1^{}\|,\label{armijot1}\\ 
					 \gamma\|\nabla_{x_1}f_1(x_1^k,x_2^k+td_2^{})\|\cdot\|g_1^{k}\| \leq \|d_1^{}\|\cdot\|g_1^k\|,\label{armijod1}\\
					 f_2(x^k+td) \leq f_2(x_1^k+td_1^{},x_2^k) + \alpha t\nabla_{x_2}f_2(x_1^k+td_1^{},x_2^k)^Td_2^{} , \label{armijo2} \\ 
					 \nabla_{x_2}f_2(x_1^k+td_1^{},x_2^k)^Td_2^{} \leq -\theta\|\nabla_{x_2}f_2(x_1^k+td_1^{},x_2^k)\|\cdot\|d_2^{}\|,\label{armijot2}\\ 
					\gamma\|\nabla_{x_2}f_2(x_1^k+td_1^{},x_2^k)\|\cdot\|g_2^{k}\| \leq \|d_2^{}\|\cdot\|g_2^k\|.\label{armijod2} 
				\end{align}  
			If one of the inequalities (\ref{armijo1})-(\ref{armijod2}) do not hold, set $t:=\frac{t}{2}$ and go to Step~2.

			\item[Step 6.] Set $t^k:=t,$ $M_1^{k}:=M_1^{},$ $M_2^{k}:=M_2^{},$ $d_1^{k}:=d_1^{},$ $d_2^{k}:=d_2^{}$ and update $x_1^{k+1}:=x_1^k+t^kd_1^{k}$, $x_2^{k+1}:=x_2^k+t^kd_2^{k}$, and $k:=k+1$. Go to Step 1.
	\end{description}
\end{breakablealgorithm}

Notice that for finding the direction $(d_1^{},d_2^{})$, we use the linear system \eqref{NewtonsistKKT2}, which for $i=1,2$ replaces the mixed Hessian matrix $t\nabla^2_{x_ix_{\neg i}}f_i(x_1^k,x_2^k)$, used in the previous discussion when deducing \eqref{NewtonsistKKT}, with zeroes whenever $g_i^k=0$ and $t\leq\tau$. This is done in order to force $d_i^{}$ to be zero if stationarity is already reached for player $i$ and the step size is sufficiently small. As a consequence, a standard Newton step is set for the other player. The reason for doing this is that in these situations it may not be possible to obtain a function decrease for player $i$. We also replaced the aforementioned condition \eqref{danterior}, which controls the ratio between gradient and direction, by \eqref{armijod1}. Both conditions are the same when $g_1^k \neq 0$, but using \eqref{armijod1} ensures that when stationarity is reached for the first player, inequality \eqref{armijod1} is automatically satisfied. A similar situation occurs for the second player with respect to condition \eqref{armijod2}.

On the other hand, taking $d_i^{k}=0$ may not be the most adequate strategy, even if $g_i^k=0$, since stationarity for player $i$ may be lost after a step is made for the other player. Thus, this is done only as a safeguarding procedure after the stepsize is at least as small as a threshold value $\tau \in (0,1]$. Notice that when $\tau:=1$, the safeguarding procedure is always activated whenever $g_i^k=0$, while when $\tau<1$, it is activated only after rejecting the Newtonian step at least once. This is done in order to try using the Newtonian direction as much as possible due to its fast local properties.

\section{Well Definiteness of the Algorithm}\label{boadef}

This section is devoted to proving well-definiteness of Algorithm \ref{Algorithm 1}. We begin by showing that when the step-size $t$ is sufficiently small the matrix of the linear system \eqref{NewtonsistKKT2} is non-singular and that its solutions can be bounded by the KKT residual.
 Before stating these results, we establish some notation that shall be used in the sequel. We denote
\begin{eqnarray*}
&  g(x):= \left[\begin{array}{c} \nabla_{x_1}f_1(x)\\ \nabla_{x_2}f_2(x) \end{array}\right], \quad g^k := \left[\begin{array}{c} g_1^k \\ g_2^k\end{array}\right] = g(x^k), \quad {\cal H}_{t,k} := \left[ \begin{array}{cc}
H_1^k & tM_1 \\ tM_2&H_2^k \end{array}\right], \\
&  d^{k} := \left[\begin{array}{c} d_1^{k} \\ d_2^{k}\end{array}\right], \quad d^{} := \left[\begin{array}{c} d_1^{} \\ d_2^{}\end{array}\right], \quad \text{ and } \quad x^{k} := \left[\begin{array}{c} x_1^{k} \\ x_2^{k}\end{array}\right].
\end{eqnarray*}

We also denote $\hat{M}^k = \max\{\|\nabla^2_{x_1x_2} f_1 (x^k)\|,\|\nabla^2_{x_2x_1} f_2 (x^k)\|\}$ where $\frac{1}{\hat{M}^k}:=+\infty$ when $\hat{M}^k=0$. The first result of  this section establishes bounds for the matrices ${\cal H}_{t,k}$ and ${\cal H}^{-1}_{t,k}$ when $t$ is sufficiently small. The following assumption will be needed:

\begin{assump}\label{Assump:Huniforposdeflim} For all $k \in \mathbb{N}$ the matrices $H^k_i$ are symmetric positive definite and there exist $\lambda_{max} \geq \lambda_{min}>0$ such that for all $i \in \{1,2\}$ and $k \in \mathbb{N}$ the eigenvalues of $H_i^k$ lie in the interval $[\lambda_{min},\lambda_{max}] $. As a consequence, for all $x_i\in \R^{n_i}$, we have 
\begin{equation*}\begin{split}\lambda_{min} \leq \|H_i^k\| \leq \lambda_{max},\quad
\frac{1}{\lambda_{max}} \leq \|(H_i^k)^{-1}\| \leq \frac{1}{\lambda_{min}}, \quad \text{ and}\\
\lambda_{min} \|x_i\|^2\leq x_i^TH_i^kx_i  \leq \lambda_{max} \|x_i\|^2,\quad i \in \{1,2\}.\end{split}\end{equation*}
\end{assump}
Assumption \ref{Assump:Huniforposdeflim} is standard in optimization and it can be guaranteed when choosing the matrices $H_i^k$. For instance, it is satisfied when the identity matrix is chosen in all iterations. This choice would be associated with the gradient method in a certain sense. Other standard choice would be considering $H_i^k=\nabla^2_{x_ix_i}f_i(x^k)$ and modifying it if necessary by altering its diagonal or by using a modified Cholesky factorization, as is done in optimization methods \cite{modcholesky}. This strategy favors the pure Newtonian direction \eqref{std} whenever possible.

\begin{lemma}\label{lemah2}
Suppose that Assumption \ref{Assump:Huniforposdeflim} holds.
For every $k \in \mathbb{N}$ and $t \in (0,1]$, the matrix ${\cal H}_{t,k}$ satisfies
\begin{equation}\label{boundsuph}
\|{\cal H}_{t,k}\| \leq \sqrt{\lambda_{max}^2 + 4\lambda_{max}\hat{M}^k + (\hat{M}^k)^2}.
\end{equation} 
Moreover, if $t \leq \frac{\lambda_{min}^2}{8\lambda_{max}\hat{M}^k}$, then ${\cal H}_{t,k}$ is non-singular and it holds
	$\|{\cal H}_{t,k}^{-1}\| \leq \frac{\sqrt{2}}{\lambda_{min}}. $
\end{lemma}
{\bf Proof:}

By the definition of ${\cal H}_{t,k}$, for a vector $z=(z_1,z_2)\in \mathbb{R}^{n_1} \times \mathbb{R}^{n_2}$, we have that
\[{\cal H}_{t,k}z = \left[\begin{array}{cc} H_1^k & tM_1\\ tM_2 & H_2^k\end{array}\right]\cdot \left[\begin{array}{c}z_1\\ z_2\end{array}\right] = \left[\begin{array}{c}H_1^kz_1 + tM_1z_2 \\ tM_2z_1 + H_2^kz_2\end{array}\right].\]

Therefore, we can write
\begin{align}
 \|{\cal H}_{t,k}z\|^2& = (H_1^kz_1 + tM_1z_2)^T(H_1^kz_1 + tM_1z_2) + (tM_2z_1 + H_2^kz_2)^T(tM_2z_1 + H_2^kz_2) \nonumber \\
& = \|H_1^kz_1\|^2 +2 tz_1^TH_1^kM_1z_2 + \|tM_1z_2\|^2 + \|H_2^kz_2\|^2  +2 tz_2^TH_2^kM_2z_1  \nonumber \\
& + \|tM_2z_1\|^2. \label{aproveitar}
\end{align}
Now, since $\|M_i\|\leq \hat{M}^k$, we have that
\begin{align}
& \|{\cal H}_{t,k}z\|^2 \leq \|H_1^kz_1\|^2 + 2t\|H_1^k\|\cdot \hat{M}^k\|z_1\|\cdot\|z_2\| \nonumber \\
&  + t^2(\hat{M}^k)^2\|z_2\|^2 + \|H_2^kz_2\|^2 + 2t\|H_2^k\|\cdot \hat{M}^k\|z_1\|\cdot\|z_2\| \nonumber \\
& + t^2(\hat{M}^k)^2 \|z_1\|^2. \nonumber
\end{align}

Using Assumption \ref{Assump:Huniforposdeflim}, and the fact that $t\leq 1$, we have that
\begin{align*}
 & \|{\cal H}_{t,k}z\|^2 \leq \lambda_{max}^2\|z_1\|^2 + 2t\lambda_{max}\hat{M}^k\|z_1\|\cdot\|z_2\| + t^2(\hat{M}^k)^2\|z_2\|^2 \nonumber \\
& + \lambda_{max}^2\|z_2\|^2 + 2t\lambda_{max}\hat{M}^k\|z_1\|\cdot\|z_2\| + t^2(\hat{M}^k)^2 \|z_1\|^2 \nonumber \\
& \leq \lambda_{max}^2\left(\|z_1\|^2 + \|z_2\|^2\right) + 4\lambda_{max}\hat{M}^k\|z_1\|\cdot\|z_2\| + (\hat{M}^k)^2\left(\|z_1\|^2 + \|z_2\|^2\right) \nonumber \\
& = \left(\lambda_{max}^2 + (\hat{M}^k)^2\right)\left(\|z_1\|^2 + \|z_2\|^2\right) + 4\lambda_{max}\hat{M}^k\|z_1\|\cdot\|z_2\|.
\end{align*}

Hence, observing that $\|z_i\| \leq \|z\|$ for $i=1,2$, we have that
\[\|{\cal H}_{t,k}z\|^2 \leq \left(\lambda_{max}^2 + 4\lambda_{max}\hat{M}^k + (\hat{M}^k)^2\right)\|z\|^2,\]
and thus \eqref{boundsuph} holds. To prove the bound on $\|{\cal H}^{-1}_{t,k}\|$, recalling \eqref{aproveitar}, we have that
\begin{align}
\|{\cal H}_{t,k}z\|^2 &=  \|H_1^kz_1\|^2 +2 tz_1^TH_1^kM_1z_2 + \|tM_1z_2\|^2 + \|H_2^kz_2\|^2  +2 tz_2^TH_2^kM_2z_1 \nonumber \\
& + \|tM_2z_1\|^2 \geq \lambda_{min}^2\|z_1\|^2 - 2t\lambda_{max}\|M_1\|\cdot\|z_1\|\cdot\|z_2\| + \|tM_1z_2\|^2  \nonumber \\
& + \lambda_{min}^2\|z_2\|^2 - 2t\lambda_{max}\|M_2\|\cdot\|z_1\|\cdot\|z_2\| + \|tM_2z_1\|^2.\nonumber
\end{align}
Since $\|tM_iz_{\neg i}\|^2\geq0$, $\|M_i\|\leq \hat{M}^k$, $\|z\|^2=\|z_1\|^2+\|z_2\|^2$, and $\|z_i\|\leq \|z\|$, we have
\begin{align}
 \|{\cal H}_{t,k}z\|^2 
& \geq \left(\lambda_{min}^2 - 4t\lambda_{max}\hat{M}^k\right)\|z\|^2. \nonumber
\end{align}
Thus, if either $\hat{M}^k=0$ or if $t \leq \frac{\lambda_{min}^2}{8\lambda_{max}\hat{M}^k}$, we have
\begin{equation}\label{otherbh}
\|{\cal H}_{t,k}z\|^2 \geq  \frac{\lambda_{min}^2}{2}\|z\|^2,
\end{equation} 
so ${\cal H}_{t,k}$ is non-singular. Finally, since  ${\cal H}_{t,k}$ is invertible, for every vector $w \in \mathbb{R}^{n_1+n_2}$ we can write $z = {\cal H}_{t,k}^{-1}w$, so \eqref{otherbh} ensures that
\begin{equation*}
\|{\cal H}_{t,k}^{-1}w\|^2 = \|z\|^2 \leq \frac{2}{\lambda_{min}^2}\|{\cal H}_{t,k}z\|^2 = \frac{2}{\lambda_{min}^2}\|{\cal H}_{t,k} {\cal H}_{t,k}^{-1}w\|^2 =  \frac{2}{\lambda_{min}^2} \|w\|^2,
\end{equation*}
implying that $\|{\cal H}_{t,k}^{-1}\|\leq\frac{\sqrt{2}}{\lambda_{min}}$.
\halmos

Next we prove a technical result that shall be used later.

\begin{lemma}\label{lemah3}
	Suppose that Assumption \ref{Assump:Huniforposdeflim} holds. For every $k \in \mathbb{N}$, let  $C_k := \frac{\sqrt{2}}{\lambda_{min}}(\|g_1^k\|+ \|g_2^k\|)$ and $d$ be as in Step 4 of Algorithm \ref{Algorithm 1}. If $t \leq \frac{\lambda_{min}^2}{8\lambda_{max}\hat{M}^k}$, then $\|d^{}\| \leq C_k$.
\end{lemma}
{\bf Proof:} By the system equation \eqref{NewtonsistKKT2} and using Lemma \ref{lemah2}, we have that 
\begin{equation*}\label{ck}
\|d^{}\| = \| {\cal H}^{-1}_{t,k} g^k\| \leq \|{\cal H}^{-1}_{t,k}\|\cdot (\|g_1^k\|+ \|g_2^k\|) \leq \frac{\sqrt{2}}{\lambda_{min}}(\|g_1^k\|+ \|g_2^k\|).
\end{equation*}
\halmos
The next result shows that for each component $i=1,2$ and each step $k$, the gradients $g_i^k$ are comparable with $d_i^{}$ provided that $t$ is small enough.

\begin{lemma}\label{lemadordemg}
	Suppose that Assumption \ref{Assump:Huniforposdeflim} holds. For every $k \in \mathbb{N}$ and $i=1,2$, if $t$ satisfies
	\begin{equation}\label{boundtlema33}
	t \leq \min\left\{ \frac{\lambda_{min}^2}{8\lambda_{max}\hat{M}^k}, \frac{\|g_i^k\|}{2\hat{M}^kC_k} \right\},
	\end{equation} where $C_k$ is as in Lemma \ref{lemah3}, then 
	\begin{equation}\label{duaslema3}
	\frac{1}{2}\|g_i^k\| \leq \|H_i^kd^{}_i\| \leq \frac{3}{2}\|g_i^k\|.
	\end{equation}
\end{lemma}
{\bf Proof:} We shall prove the result for $i=1$, being the result analogous for $i=2$. Using the triangle inequality in the  first equation of system \eqref{NewtonsistKKT2} we have that
\begin{equation}\label{firstlema3}
 \|g_1^k\| - \|tM_1d_2^{} \| \leq \|H_1^kd_1^{}\| \leq \|tM_1d_2^{}\|+\| g_1^k\| .
\end{equation}
Thus, if $\hat{M}^k=0$, then $M_1=0$ and we have nothing to prove. Otherwise, by Lemma \ref{lemah3}, \eqref{firstlema3} implies that
\begin{equation*}
 \|g_1^k\| - t\hat{M}^k C_{k} \leq \|H_1^kd_1^{}\|  \leq t\hat{M}^k C_{k} + \|g_1^k\|.
\end{equation*}
So, from \eqref{boundtlema33}, we conclude that \eqref{duaslema3} holds.
\halmos

As previously mentioned, by the way Algorithm \ref{Algorithm 1} is built, in cases where one of the gradients is zero, say $\|g_2^k\|=0$, and the step-size $t$ has sufficiently decreased ($t\leq\tau$), it turns out that the direction for the second player is null while for the first player the direction is computed as a standard Newton direction for optimization with Armijo line-search, which is known to be well-defined if the objective functions have locally Lipschitz continuous gradients with respect to the decision variables $x_i$ for every fixed $x_{\neg i}$, $i = 1,2$. In order to be more precise, let us consider the following hypothesis:
\begin{assump} \label{Assump:gradlocLipz}
    For each $i\in\{1,2\}$ and $x = (x_i,x_{\neg i}) \in \R^{n_1+n_2}$, there exist $\delta_x,L_x>0$ such that  
    $$\|\nabla_{x_i} f_i(y, x_{\neg i})-\nabla_{x_i} f_i(z, x_{\neg i})\| \leq L_x \|y-z\|, \quad \forall y,z \in B(x_i,\delta_x).$$ 
\end{assump}
As a consequence, we have that
	\begin{equation}\label{boundquadraticLips}
	 f_i(x_i+td_i,x_{\neg i}) \leq f_i(x) + t\nabla_{x_i} f_i(x)^Td_i + \frac{t^2L_x}{2}\|d_i\|^2
    \end{equation}
		for every  $x = (x_i, x_{\neg i})\in \mathbb{R}^{n_1+n_2}, d_i \in \mathbb{R}^{n_i}$, and $t>0$ such that $\|td_i\| \leq \delta_x$, $i\in\{1,2\}$.

For $x=x^k$ generated by Algorithm \ref{Algorithm 1}, we denote $\delta_x=\delta_k$ and $L_x=L_k.$

Note that, when $g_1^k=0$ and $t\leq \tau$, then by \eqref{separar} and \eqref{NewtonsistKKT2} we have $d_1=0$. In this case,  conditions \eqref{armijo1}-\eqref{armijod1} are trivially satisfied. Furthermore, $d^k_2$ is the standard Newtonian direction for $f_2(x_1^k,x_2)$ at $x_2=x_2^k$. So, under Assumptions \ref{Assump:Huniforposdeflim} and \ref{Assump:gradlocLipz}, conditions \eqref{armijo2}-\eqref{armijod2}  hold for all  \begin{equation}\label{tafstado0otim}
    t~\leq~ \min\left\{\frac{2(1-\alpha)\lambda_{min}^2}{L_k\lambda_{max}},\frac{\lambda_{min}\delta_k}{\max\{\|g_2\|,1\}}\right\},
\end{equation}  whenever $\theta = \frac{\lambda_{min}}{\lambda_{max}}$ and $\gamma=\frac{1}{\lambda_{max}}$ (see for instance the discussion in \cite{bertsekas}, pages 29-36, and Lemma 2.20 in \cite{solodov}). Analogous situation shows that the iteration is also well defined if $g_2=0$. Thus, we focus now on the case where $\|g_i^k\|\neq 0, i=1,2$, and we show that the inequalities (\ref{armijo1})-(\ref{armijod2}) are satisfied for sufficiently small $t$, meaning that the step-size does not decrease indefinitely.

We begin with the gradient/direction ratio inequalities \eqref{armijod1} and \eqref{armijod2}. However, we show a stronger version bounding the directions by the gradient on the points $(x_i^k, x_{\neg i}^k + td_{\neg i}^{})$ from both above and below. In order to achieve this, we use the continuity of the second derivatives of $f_i$, $i =1,2$, which guarantees that for every $k \in \N$, 
there exists $\bar{\delta}_k>0$ such that
\begin{equation}\label{conthessmista}
  \|\nabla^2_{x_ix_{\neg i}} f_i(x)\| \leq  \hat{M}^k +1, \quad \forall x \in B(x^k, \bar{\delta}_k).
\end{equation}

\begin{lemma}\label{infereta}
Suppose that Assumption \ref{Assump:Huniforposdeflim} holds.
For every $k \in \mathbb{N}$ and $i=1,2$, if $t$ satisfies
\begin{equation}\label{boundtlema1}
t \leq \min\left\{\frac{\lambda_{min}^2}{8\lambda_{max}{\hat{M}^k}}, \frac{\|g_i^k\|}{4(2\hat{M}^k+1)C_k},  \frac{\bar{\delta}_k}{C_k}\right\}, 
\end{equation}
where  $C_{k}$ is as in Lemma \ref{lemah3} and $\bar{\delta}_k$ is as in \eqref{conthessmista}, then 
\begin{equation}\label{etacondition}
\frac{2}{3\lambda_{max}} \left\Vert \nabla_{x_i} f_i(x_i^k, x_{\neg i}^k + td_{\neg i}^{})\right\Vert \leq \|d^{}_i\| \leq \frac{2}{\lambda_{min}} \left\Vert \nabla_{x_i} f_i(x_i^k, x_{\neg i}^k + td_{\neg i}^{})\right\Vert.
\end{equation}

\end{lemma}
{\bf Proof:} We prove the result only for $i=1$, as the proof for $i=2$ is similar. By~\eqref{NewtonsistKKT2}, we have that
\begin{equation}\label{eqsist1}
-g_1^k = H_1^kd_1^{} + tM_1d_2^{}.
\end{equation}
Also, we have that
\begin{equation}\label{eq2}
\nabla_{x_1}f_1(x_1^k,x_2^k + td_2^{}) = g_1^k + t\left[ \displaystyle \int_0^1 \nabla_{x_1x_2}^2f_1(x_1^k,x_2^k + \xi td_2^{})d\xi \right] d_2^{}.
\end{equation}
Combining \eqref{eqsist1} and \eqref{eq2} we obtain
\begin{equation}\label{toreorder}
\nabla_{x_1}f_1(x_1^k,x_2^k + td_2^{}) = -H_1^kd_1^{} + t\left[ -M_1 + \displaystyle \int_0^1 \nabla_{x_1x_2}^2f_1(x_1^k,x_2^k + \xi td_2^{})d\xi \right] d_2^{}.
\end{equation}
By the triangular inequality, we have both inequalities: 
\begin{align}
 & \|\nabla_{x_1}f_1(x_1^k,x_2^k + td_2^{})\| \geq  \Bigg| \|H_1^kd_1^{}\| - t\bigg\lVert\left[ -M_1 + \displaystyle \int_0^1 \nabla_{x_1x_2}^2f_1(x_1^k,x_2^k + \xi td_2^{})d\xi \right]d_2^{} \bigg\rVert \Bigg| \label{maineq} \\
&\|\nabla_{x_1}f_1(x_1^k,x_2^k + td_2^{})\| \leq  \|H_1^kd_1^{}\|+t\bigg\lVert\left[ -M_1 + \displaystyle \int_0^1 \nabla_{x_1x_2}^2f_1(x_1^k,x_2^k + \xi td_2^{})d\xi \right]d_2^{} \bigg\rVert. \label{maineq2}
\end{align}
Next, we prove that for $t$ sufficiently small it holds
\begin{equation}\label{middesigf}
\frac{1}{2}\|H_1^kd_1^{}\| \leq \|\nabla_{x_1}f_1(x_1^k,x_2^k + td_2^{})\| \leq \frac{3}{2}\|H_1^kd_1^{}\|.
\end{equation}

Let us suppose that $\left[ -M_1 + \displaystyle \int_0^1 \nabla_{x_1x_2}^2f_1(x_1^k,x_2^k + \xi td_2^{})d\xi \right]d_2^{} \neq 0$, since otherwise the result is trivial from \eqref{toreorder}. Hence, if we select $t$ satisfying 
\begin{equation}\label{tsatisfaz1}
t \leq \frac{\|H_1^kd_1^{}\|}{2 \bigg\lVert\left[ -M_1 + \displaystyle \int_0^1 \nabla_{x_1x_2}^2f_1(x_1^k,x_2^k + \xi td_2^{})d\xi \right]d_2^{} \bigg\rVert},
\end{equation}
we see from \eqref{maineq} that the left-hand side inequality in \eqref{middesigf} holds and, by \eqref{maineq2}, we obtain the right-hand side inequality. Now, using that $\lambda_{min}\|d_1^{}\| \leq \|H_1^kd_1^{}\| \leq \lambda_{max}\|d_1^{}\|$, \eqref{middesigf} implies that
\[\frac{\lambda_{min}}{2}\|d_1^{}\| \leq \|\nabla_{x_1}f_1(x_1^k,x_2^k + td_2^{})\| \leq \frac{3\lambda_{max}}{2}\|d_1^{}\|,\]
which implies the desired bounds in \eqref{etacondition}.

Finally, let us show that $t$ satisfying \eqref{boundtlema1} implies the bound \eqref{tsatisfaz1}. Since \eqref{boundtlema1} and Lemma \ref{lemah3} imply that $\|d_2\| \leq \|d\| \leq C_k$, we have that $t \|d_2\| \leq \bar{\delta}_k$. By \eqref{conthessmista} , we see that 
\begin{align*}
&\bigg\lVert\left[ -M_1 + \displaystyle \int_0^1 \nabla_{x_1x_2}^2f_1(x_1^k,x_2^k + \xi td_2^{})d\xi \right]d_2^{} \bigg\rVert\\
& \leq \left[\| M_1 \| + \displaystyle \int_0^1 \lVert \nabla_{x_1x_2}^2f_1(x_1^k,x_2^k + \xi td_2^{}) \rVert d\xi \right] \cdot \|d_2^{}\| \leq  (2\hat{M}^k+1)C_k.
\end{align*}
 Therefore, in order to attain \eqref{tsatisfaz1}, it suffices to select
\begin{equation}\label{tboundlema}
t \leq \frac{\|H_1^kd_1^{}\|}{2(2\hat{M}^k+1)C_k}.
\end{equation}
Now, notice that inequalities \eqref{duaslema3} guaranteed by Lemma \ref{lemadordemg} ensure \eqref{tboundlema} holds when $t$ satisfies \eqref{boundtlema1}, which concludes the proof.
\halmos

The next result ensures the angle inequalities \eqref{armijot1} and \eqref{armijot2} hold for $t$ sufficiently small, provided that we require some control on the linear approximation of the gradients, as is stated in the assumption below.

\begin{assump}\label{Asssump:RestoGrad}
For every $x^k \in \R^{n_1+n_2}$ generated by Algorithm \ref{Algorithm 1}, there exist $\hat{\delta}_k>0,C_R^k>0$, and functions $r_i^k:\R^{n_{\neg i+1}} \to \R$, $i =1,2$, such that  for all $d_{\neg i} \in \mathbb{R}^{n_{\neg i}}$ and for $t>0$  satisfying $\|t d_{\neg i}\| \leq \hat{\delta}_k$ it holds 
\begin{equation}\label{restotaylor1grad}
	\nabla_{x_i} f_i(x_i^k,x_{\neg i}^k+td_{\neg i}) = \nabla_{x_i} f_i(x^k)+t\nabla_{x_ix_{\neg i}}^2f_i(x^k)d_{\neg i}+r_i^k(t,d_{\neg i}),
\end{equation}
and
	$\| r_i^k(t,d_{\neg i})\| \leq C_R^k t^2 \|d_{\neg i}\|^2$.
\end{assump}
We highlight that assumption \ref{Assump:gradlocLipz} and \ref{Asssump:RestoGrad} are simple requirements on the smoothness of the high-order derivatives of the functions $f_1$ and $f_2$.

\begin{lemma}\label{prop_angle} Suppose that Assumptions \ref{Assump:Huniforposdeflim} and \ref{Asssump:RestoGrad} hold. For every $k \in \mathbb{N}$ and $i=1,2$, if $g_i^k\neq 0$ and $t$ satisfies
\begin{equation}\label{t1propv}
t \leq \min\left\{\frac{\lambda_{min}^2}{8\lambda_{max}{\hat{M}^k}}, \frac{\|g_i^k\|}{4(2\hat{M}^k+1)C_k},\displaystyle \sqrt{\frac{\lambda_{min}\|g_i^{k}\|}{8\lambda_{max}C_R^k(C_k)^2}} ,  \frac{\bar{\delta}_k}{C_k},\frac{\hat{\delta}_k}{C_k} \right\},
\end{equation}
where $C_{k}$ is given as in Lemma \ref{lemah3} and $\bar{\delta}_k$, $\hat{\delta}_k$ are given as in \eqref{conthessmista} and Assumption \ref{Asssump:RestoGrad}, respectively. Then 
\begin{equation}\label{angle_condition}
\nabla_{x_i} f_i(x^k_i,x^{k}_{\neg i}+td_{\neg i}^{})^Td_i^{} \leq -\frac{\lambda_{min}}{4\lambda_{max}} \|\nabla f_i(x^k_i,x^{k}_{\neg i}+td_{\neg i}^{})\|\cdot\|d_i^{}\|.
\end{equation}
\end{lemma}
{\bf Proof:} Once again, we work the proof of the main claim for the first component. For the second one, the proof follows similarly. Since we are considering $g_1^k \neq 0$, we have that $M_1 = \nabla_{x_1x_2}^2f_1(x^k)$. Hence, Algorithm \ref{Algorithm 1} yields
\begin{equation}\label{eqmainstep}
H_1^k d_1^{} = - g_1^k - t\nabla_{x_1x_2}^2f_1(x^k)d_2^{}.
\end{equation}

Since \eqref{t1propv} implies that $\|t d_2\| \leq \hat{\delta}_k$, by Assumption \ref{Asssump:RestoGrad}, we have that
\begin{equation}\label{firstorder}
\nabla_{x_1}f_1(x_1^k,x_2^{k}+td_2^{}) = g_1^k + t\nabla^2_{x_1x_2}f_1(x^k)d_2^{} + r_1^k(t,d_2^{}), 
\end{equation}
where the remainder $ r_1^k(t,d_2^{})$ is such that
$\| r_1^k(t,d_2^{})\| \leq C_R^kt^2\|d_2^{}\|^2$.

Therefore, adding \eqref{firstorder} to \eqref{eqmainstep}, we have
\begin{equation}\label{gives}
 H_1^kd_1^{} = -\nabla_{x_1} f_1(x_1^k,x_2^{k}+td_2^{}) +  r_1^k(t,d_2^{}). 
\end{equation}
Hence, since $\|H_1^kd_1^{}\|\leq \lambda_{max}\|d_1^{}\|$,
 using the reverse triangule inequality, 
we obtain
\begin{equation}\label{guardarprop2}
\|d_1^{}\| \geq \frac{1}{\lambda_{max}}\Big| \| r_1^k(t,d_2^{})\| - \|\nabla_{x_1}f_1(x_1^k,x_2^{k}+td_2^{})\| \Big|.
\end{equation} 
On the other hand, multiplying \eqref{gives} by $d_1^{T}$ 
and using the Cauchy-Schwarz inequality, Assumption \ref{Assump:Huniforposdeflim} and Lemma \ref{lemah3}, we have that
\begin{align}\label{citadod}
 d_1^T\nabla_{x_1}f_1(x_1^k,x_2^{k}+td_2^{}) &\leq -d_1^TH_1^kd_1^{} + \| r_1^k(t,d_2^{})\| \cdot \|d_1^{}\| \nonumber \\
 & \leq -\lambda_{min}\|d_1^{}\|^2 + C_R^kC_k^2t^2\|d_1^{}\|.
\end{align}

Therefore, if we choose $t$ in \eqref{citadod} satisfying
$t \leq \displaystyle \sqrt{\frac{\lambda_{min}\|d_1^{}\|}{2C_R^kC_k^2}},$
we obtain
\begin{equation}\label{sobre2}
d_1^T\nabla_{x_1}f_1(x_1^k,x_2^{k}+td_2^{}) \leq \frac{-\lambda_{min}\|d_1^{}\|^2}{2}.
\end{equation}

Thus, we can combine \eqref{guardarprop2} and \eqref{sobre2}, to write
\begin{equation}
d_1^T\nabla_{x_1}f_1(x_1^k,x_2^k+td_2^{}) \leq -\frac{\lambda_{min}}{2\lambda_{max}}\|d_1^{}\|\Big| \|\nabla_{x_1}f_1(x_1^k,x_2^k+td_2^{})\| - \| r_1^k(t,d_2^{})\| \Big|. \label{ineqmodulo}
\end{equation}
Next, to bound the right-hand side of \eqref{ineqmodulo}, for $t$ satisfying \eqref{t1propv} we have by \eqref{duaslema3} that
\begin{align}
 t^2 \leq \frac{\lambda_{min}\|g_1^{k}\|}{8\lambda_{max}C_R^kC_k^2} \leq  \frac{\lambda_{min}\|H_1^kd_1^{}\|}{4\lambda_{max}C_R^kC_k^2} \leq \frac{\lambda_{min}\|d_1^{}\|}{4C_R^kC_k^2}. \label{boundmediot}
\end{align}
Hence, using \eqref{boundmediot} in the bound of $r_1^k(t,d_2)$ given in Assumption \ref{Asssump:RestoGrad}, we have that
\[\| r_1^k(t,d_2^{})\| \leq C_R^kt^2\|d_2^{}\|^2 \leq C_R^kC_k^2t^2 \leq \frac{\lambda_{min}\|d_1^{}\|}{4} \leq \frac{\|\nabla_{x_1}f_1(x_1^k,x_2^k+td_2^{})\|}{2},\]
where we used Lemma \ref{infereta} in the last inequality. Therefore, we have that
\begin{equation*}
	\Big| \|\nabla_{x_1}f_1(x_1^k,x_2^k+td_2^{})\| -\| r_1^k(t,d_2^{})\| \Big| \geq \frac{\|\nabla_{x_1}f_1(x_1^k,x_2^k+td_2^{})\|}{2}.
\end{equation*}
Then, \eqref{ineqmodulo} gives us the desired inequality
\[d_1^T\nabla_{x_1}f_1(x_1^k,x_2^k+td_2^{}) \leq -\frac{\lambda_{min}}{4\lambda_{max}}\|d_1^{}\|\cdot \|\nabla_{x_1}f_1(x_1^k,x_2^k+td_2^{})\|.\]
\halmos

Now we are ready to finish establishing the well-definiteness of Algorithm \ref{Algorithm 1}, where we prove the remaining relations \eqref{armijo1} and \eqref{armijo2}.

\begin{proposition}\label{welldefinedblue}
Suppose that Assumptions \ref{Assump:Huniforposdeflim}, \ref{Assump:gradlocLipz}, and \ref{Asssump:RestoGrad} hold. For $i=1,2$, suppose that $g_i^k\neq 0$ for some $k \in \mathbb{N}$, and assume that the step-size $t$ satisfies
\begin{equation}\label{boundtblue}
t \leq \min\left\{\frac{\lambda_{min}^2}{8\lambda_{max}{\hat{M}^k}}, \frac{\|g_i^k\|}{4(2\hat{M}^k+1)C_k},  \displaystyle \sqrt{\frac{\lambda_{min}\|g_i^{k}\|}{8\lambda_{max}C_R^kC_k^2}},\frac{(1-\alpha)\lambda_{min}^2}{4 L_k\lambda_{max}}, \frac{\delta_k}{C_k},\frac{\bar{\delta}_k}{C_k}, \frac{\hat{\delta}_k}{C_k} \right\},  
\end{equation}
where $\delta_k$, $\bar{\delta}_k$ and $\hat{\delta}_k$ are given as in Assumption \ref{Assump:gradlocLipz}, \eqref{conthessmista} and Assumption \ref{Asssump:RestoGrad}, respectively. 
Then, if
\begin{equation}\label{valoresgammatheta}
\gamma \leq \frac{2}{3\lambda_{max}} \quad \text{and} \quad \theta \leq \frac{\lambda_{min}}{4\lambda_{max}},
\end{equation}
all inequalities (\ref{armijo1})-(\ref{armijod2}) are satisfied. As a consequence, with the choice of parameters $\gamma$ and $\theta$ satisfying \eqref{valoresgammatheta} and any $\alpha \in(0,1)$, 
Algorithm \ref{Algorithm 1} is well-defined.
\end{proposition}

{\bf Proof:}
By Lemmas \ref{infereta} and \ref{prop_angle}, it remains to prove \eqref{armijo1}, since \eqref{armijo2} is analogous.
Let us, fix $k \in \mathbb{N}$ and $i=1$. Since \eqref{boundtblue} implies that $\|td_2\| \leq \delta_k$, by Assumption \ref{Assump:gradlocLipz}, we have that
\begin{equation*}
f_1(x_1^{k}+td_1^{},x_2^{k}+td_2^{}) - f_1(x_1^k,x_2^{k}+td_2^{}) \leq t\nabla_{x_1}f_1(x_1^k,x_2^{k}+td_2^{})^Td_1^{} + \frac{t^2L_k}{2}\|d_1^{}\|^2.
\end{equation*}

If $d_1^{} = 0$, then the result is immediate. If $d_1^{} \neq 0$ and $t>0$, we have that
\begin{align*}
& t\nabla_{x_1}f_1(x_1^k,x_2^k+t d_2^{})^Td_1^{} + \frac{t^2L_k}{2}\|d_1^{}\|^2 \leq \alpha t \nabla_{x_1}f_1(x_1^k,x_2^k+t d_2^{})^Td_1^{} \nonumber \\
& \Leftrightarrow t \leq \frac{2(1 - \alpha)\left[ -\nabla_{x_1}f_1(x_1^k,x_2^k+t d_2^{})^Td_1^{} \right]}{L_k\|d_1^{}\|^2},\label{semteta}
\end{align*}
which can be guaranteed by 
$t \leq \frac{2(1-\alpha)\frac{\lambda_{min}}{4\lambda_{max}}\|\nabla_{x_1}f_1(x_1^k,x_2^k+t d_2^{})\|}{L_k\|d_1^{}\|}$
due to \eqref{angle_condition}. Also, by Lemma \ref{infereta}, we have that $\frac{\|\nabla_{x_1}f_1(x_1^k,x_2^k+t d_2^{})\|}{\|d_1^{}\|} \geq \frac{\lambda_{min}}{2}$, and thus, it suffices to select $t$ satisfying \eqref{boundtblue} to ensure \eqref{armijo1}.
\halmos

Notice that our well-definiteness results assume that $t$ is smaller than a positive constant which depends on the gradients $g_i^k$ at iteration $k$. In other words, once $g_i^k \neq 0$ for $i=1,2$, there exists a fixed positive value of $t$ (depending on $k$) that guarantees that the iteration will necessarily end. On the other hand, this is not guaranteed if $g_i^k = 0$ for some $i \in \{1,2\}$. This is the technical reason why we must switch to the standard Newton method for optimization whenever one player has already reached stationarity. We could not prove or find a counter-example for the well-definiteness of the method that never switches to standard Newton, that is, considering $M_i$ always equals to $\nabla^2_{x_ix_{\neg i}}f_i(x^k),$ $i=1,2,$ in place of \eqref{separar}.  With the use of this safeguarding procedure we have that any iteration of the algorithm is well defined.

We note that the Armijo-type conditions \eqref{armijo1} and \eqref{armijo2} can be guaranteed to be satisfied at a small enough $t$ without requiring local Lipschitz continuity of the gradients. On the other hand, since our algorithm uses second-order derivatives, we chose to use the Lipschitz continuity assumption (Assumption \ref{Assump:gradlocLipz}) without loss of generality. The advantage in doing so is that we can make the bounds for $t$ explicit. This is analogous to what is known in traditional optimization. On the other hand, conditions \eqref{armijot1}, \eqref{armijot2}, \eqref{armijod1}, and \eqref{armijod2}, which bound the size of the directions and their angles with the respective gradients, can be always satisfied in traditional optimization by simply controlling the eigenvalues of the matrix that approximates the Hessian of the objective function in a quasi-Newton scheme. The situation is different in our algorithm since the conditions involve the gradient of one player's objective function, but evaluated at the predicted decision of the other player, which depends on the value of $t$. Therefore, the analysis is more complex, being necessary to test conditions \eqref{armijot1}, \eqref{armijot2}, \eqref{armijod1}, and \eqref{armijod2} inside the algorithm, since they are only guaranteed to be fulfilled when $t$ is sufficiently small.

 Next, we discuss the case when the iterates all lie on a bounded set and the problem is sufficiently smooth. Due to the descent characteristic of the proposed algorithm, this is somehow expected to hold.

\begin{proposition}\label{constantesuniformes}
 Suppose that Assumptions \ref{Assump:Huniforposdeflim}, \ref{Assump:gradlocLipz}, and \ref{Asssump:RestoGrad} hold and that the sequence $\{x^k\}$ generated by Algorithm \ref{Algorithm 1} lies on a convex compact set $\Omega$. Suppose also that the third derivatives of the objective functions always exist and are continuous, then there exist universal constants $C_d$, $C_H$ and $C_R$ such that $0<C_k \leq C_d$, $\hat{M}^k<C_H$, and $C_R^k \leq C_R$, where $C_k$ is defined as in Lemma \ref{lemah3} and $C_R^k$ are the constants given by Assumption \ref{Asssump:RestoGrad}. Moreover, there exist $L>0$ such that $L_k < L$ for $L_k=L_{x^k}$ in Assumption \ref{Assump:gradlocLipz} and $\delta>0$ such that $\min\{\delta_k,\bar{\delta}_k,\hat{\delta}_k\} > \delta$ where $\delta_k$, $\bar{\delta}_k$, and $\hat{\delta}_k$ are given as in Assumption \ref{Assump:gradlocLipz}, \eqref{conthessmista}, and Assumption \ref{Asssump:RestoGrad}, respectively.
\end{proposition}
{\bf Proof:} The continuity of the first and second derivatives of $f_i(x)$, $i=1,2$, and the compactness of $\Omega$ guarantee the existence of $C_d$ and $C_H$. Note now that $\bar{\Omega}:=\cup_{x \in \Omega} B(x,1)$ is also compact. Thus, for $i=1,2$, $\|\nabla^2_{x_ix_{\neg i}}f_i(x)\|$ is uniformly continuous on $\bar{\Omega}$. Therefore, there exists $\bar{\delta} \in (0,1]$ such that \eqref{conthessmista} holds. Using the mean value inequality combined with the continuity of the second and third derivatives of $f_i(x)$, $i=1,2$, and the compactness of $\cup_{x \in \Omega} B(x,\delta)$, we obtain the existence of $C_R$ and $L$ satisfying the desired conditions.
\halmos 

The next result states we can infer the boundedness of the step-size sequence from the boundedness of both gradients from below. This guarantees that in this situation the line-search criteria will not get increasingly hard to be satisfied, forcing small steps and stalling the computation of the iterates, whenever the decisions of both players are far from stationarity. 

\begin{proposition}\label{coroexplica}
 Suppose that Assumptions \ref{Assump:Huniforposdeflim}, \ref{Assump:gradlocLipz}, and \ref{Asssump:RestoGrad} hold and that the sequence $\{x^k\}$ generated by Algorithm \ref{Algorithm 1}, with $\gamma$ and $\theta$ satisfying
\eqref{valoresgammatheta},  lies on a bounded set. Suppose also that the third derivatives of the objective functions always exist and are continuous. If there exists a constant $c>0$ and $k_0 \in \N$ such that $\|g_i^k\|> c$ for each component $i=1,2$ and all $k>k_0$, then the step-size sequence $\{t^k\}$ is bounded away from zero. 
\end{proposition}

{\bf Proof:} Let $0<\bar{t}:=\min\{t^k\}$ for $k\leq k_0$. For $k>k_0$,  using that $\|g_i^k\|> c$ and the bounds given by Proposition \ref{constantesuniformes}, Proposition \ref{welldefinedblue} ensures that if  
\[
t \leq t':=\min\left\{ \frac{\lambda_{min}^2}{8\lambda_{max}C_H},  \frac{c}{8C_HC_d+4C_d}, \displaystyle \sqrt{\frac{\lambda_{min}c}{8\lambda_{max}C_RC_d^2}}, \frac{(1-\alpha)\lambda_{min}^2 }{4\lambda_{max} L}, \frac{\delta}{C_d} \right\},
\]
then conditions (\ref{armijo1})-(\ref{armijod2}) hold. Thus, by the stepsize update rule we have that, for all $k$, 
$t^k \geq \min\left\{1,\frac{t'}{2},\bar{t}\right\}>0$.
\halmos

\section{Convergence Analysis}\label{sec3}
In this section we establish our convergence results for Algorithm \ref{Algorithm 1}. In practice, the reduction of the functional values at each iteration favors obtaining minimizers rather than other stationary points. However, as usually happens for classical descent algorithms in optimization, our convergence results are only associated with the stationarity of the limit point of the generated sequence, since we are not using more elaborated elements associated to second order conditions.

We start with our main convergence result, which  ensures, under certain assumptions, the stationarity of limit points for at least one of the players.

\begin{theorem} \label{pro:conv} Consider the sequences $\{x^k\}$ and $\{t^k\}$ generated by Algorithm \ref{Algorithm 1}. Then, the following holds:
\begin{itemize}
\item[i)] If $\lim_{\mathbb{N}'}d^k=0$ and if the mixed Hessians  are uniformly bounded on an infinite subset $\mathbb{N}'\subset \mathbb{N}$, then the sequence of gradients $\{g_1^k\}$ and $\{g_2^k\}$ converge to zero on the same subset.

\item[ii)] If Assumptions \ref{Assump:Huniforposdeflim}, \ref{Assump:gradlocLipz} and \ref{Asssump:RestoGrad} hold and the whole sequence $\{x^k\}$ converges to a point $x^*$, then either $\nabla_{x_1}f_1(x^*)=0$ or $\nabla_{x_2}f_2(x^*)=0$.

\item[iii)] Suppose that Assumptions \ref{Assump:Huniforposdeflim}, \ref{Assump:gradlocLipz}  and \ref{Asssump:RestoGrad} hold, and that the third derivatives of the objective functions always exist and are continuous. If $\{x^k\}$ converges to a point $x^*$ and if in addition the sequence of step-sizes $\{t^k\}$ is bounded away from zero, then $x^*$ is stationary for the NEP  \eqref{PJ1}, that is, $\nabla_{x_1}f_1(x^*)=0$ and $\nabla_{x_2}f_2(x^*)=0$. 
\end{itemize}
\end{theorem}
{\bf Proof:} By \eqref{NewtonsistKKT2} and \eqref{boundsuph}, we have that
\begin{equation*} 
\|g^k\| = \|{\cal H}_{t,k}d^{k}\| \leq \|{\cal H}_{t,k}\|\cdot\|d^{k}\| \leq \sqrt{\lambda_{max}^2 + 4\lambda_{max}\hat{M}^k + (\hat{M}^k)^2}\|d^{k}\|.
\end{equation*}
Since there exists an universal bound for $\hat{M}^k$ on $\mathbb{N}'$, if $\lim_{k \in \mathbb{N}'}d^k = 0$, then it holds $\lim_{k \in \mathbb{N}'}g^k = 0$, which proves item (i).
	
To prove item (iii), assume that the sequence converges; and then it lies  on a compact set. Therefore, by continuity of the mixed Hessians, there is an uniform bound on $\hat{M}^k$. Moreover, we have that
$0=\lim_{k \to \infty} x^{k+1}-x^k= \lim_{k \to \infty} t^k d^{k}.$

Since $t^k \geq \bar{t}$ for some constant $\overline{t}>0$, we have that $d^{k}$ converges to zero. Thus, by item (i) the whole gradient sequence converges to zero, implying that the limit point $x^*$ is stationary.

Finally, to prove (ii), assume that it does not hold. Thus, there exists $\epsilon>0$ such that $\|g_i^k\|\geq\epsilon$ for $i=1,2$ and $k$ large enough, say $k>k_0 \in \mathbb{N}$. Additionally, since $\{x^k\}$ is convergent, the sequence lies on a compact set and therefore $\hat{M}^k$ is uniformly bounded. Hence, by Proposition \ref{coroexplica}, we have that $\{t^k\}$ is bounded away from zero, which by item (iii) implies that $g^k \rightarrow 0$, a contradiction.
\halmos

Note that Theorem \ref{pro:conv} still allows a sequence that converges to a point which is stationary only for one of the players. This happens because the step-size may become asymptotically null when one of the gradients approaches zero. In this situation, the other player's search direction may still be large even though the displacement obtained in the iteration is small, due to how close to stationarity the other player is. It is then natural to ask whether or not it is actually possible  to show that the step-sizes remain always bounded away from zero. One way to ensure this is by guaranteeing a positive lower bound for the expression on the right-hand side of \eqref{boundtblue}. Assuming the sequence generated by the algorithm lies within a compact set $\Omega$, and that the functions $f_i(x)$, $i=1,2$, have third-order continuous derivatives, Proposition \ref{constantesuniformes} provides bounds for $\hat{M}^k$, $C_R^k$, $L_k$, $\delta_k$, $\bar{\delta}_k$, $\hat{\delta}_k$, and $C_k$. Thus, it suffices to require that $\frac{\|g_i^k\|}{C_k}$ is bounded away from zero. A sufficient condition for this is that there exist $\beta_1,\beta_2>0$ such that $\beta_1 \|\nabla_{x_1} f_1(x) \| \leq \|\nabla_{x_2} f_2(x) \| \leq \beta_2 \|\nabla_{x_1} f_1(x) \|$, for all $x \in \Omega$.

Another approach is to check the impact of the approximation used in \eqref{aproxgrad} on the algorithm's design. Using the remainder term from the approximation given in \eqref{restotaylor1grad}, the system \eqref{NewtonsistKKT2} implies that \begin{equation}\label{gives2} H_i^kd_i^{} = -\nabla_{x_i} f_i(x_i^k,x_{\neg i}^{k}+td_{\neg i}^{}) + r_i^k(t,d_{\neg i}^{}). \end{equation} To ensure the quality of $d_i$, it is natural to require that the remainder term is not significant with respect to the gradient, as done in the following theorem. 

\begin{proposition}
\label{teo:restobetatafastado0} Suppose that Assumptions \ref{Assump:Huniforposdeflim}, \ref{Assump:gradlocLipz}, and \ref{Asssump:RestoGrad} hold. Assume also that the sequence $\{x^k\}$ generated by Algorithm \ref{Algorithm 1} lies in a bounded set and that the third derivatives of the objective functions exist and are continuous. Moreover, if there exist $\delta'>0$ and $0\leq\beta <\lambda_{min}/(\lambda_{max}+\lambda_{min})$ such that
\begin{align}
    \|r_i^k(t,d_{\neg i}^{})\| \leq \beta \| \nabla_{x_i} f_i(x_i^k,x_{\neg i}+td_{\neg i})\|,\label{eq:rmenorbetagrad}\\
    0<\gamma \leq \frac{1-\beta}{\lambda_{max}} \;\; \text{ and }\;\; 0<\theta \leq  \frac{\lambda_{min}-(\lambda_{max}+\lambda_{min})\beta}{\lambda_{max}}, \label{gammathetacombeta}
\end{align}
 for all $i=1,2$, $k \in \N$, $t \in [0,\delta']$ and $|td_i| \leq \delta$, where $\delta$ is given in Proposition \ref{constantesuniformes}, then ${t^k}$ is bounded away from zero. 
\end{proposition}

{\bf Proof:} For $C_H$, $C_d$, $L$, and $\delta$ given by Proposition \ref{constantesuniformes}, let $t \in [0,t']$, where 
    \[ t' :=\min\left\{ \frac{\lambda_{min}^2}{8\lambda_{max}C_H}, \frac{2(1-\alpha)\theta \lambda_{min}}{(1+\beta)L} , \frac{\delta}{C_d}, \delta' \right\}.\]
   Since $t \leq t'$, we can combine \eqref{gives2}, \eqref{eq:rmenorbetagrad}, and the triangle inequality, obtaining
\begin{equation}
    \label{reladgradbeta}
    \frac{1-\beta}{\lambda_{max}} \| \nabla_{x_i} f_i(x_i^k,x_{\neg i}^{k}+td_{\neg i})\|
    \leq \|d_i\|
    \leq \frac{1+\beta}{\lambda_{min}}\| \nabla_{x_i} f_i(x_i^k,x_{\neg i}^{k}+td_{\neg i})\|.
\end{equation}
Using now \eqref{gives2}, \eqref{eq:rmenorbetagrad}, and the left-hand side of \eqref{reladgradbeta}, we have that
\begin{equation*}
    \label{condangbeta}
     \nabla_{x_i} f_i(x_i^k,x_{\neg i}^{k}+td_{\neg i})^Td_i
    \leq -\left( \frac{\lambda_{min}-(\lambda_{max}+\lambda_{min})\beta}{\lambda_{max}}\right) \|d_i\| \|\nabla_{x_i} f_i(x_i^k,x_{\neg i}^{k}+td_{\neg i})\|.
\end{equation*}
Thus, conditions \eqref{armijot1}, \eqref{armijod1}, \eqref{armijot2}, and \eqref{armijod2} hold for $\gamma$ and $\theta$ satisfying \eqref{gammathetacombeta}.
Moreover, using \eqref{boundquadraticLips}, \eqref{armijot1}, \eqref{armijot2}, and the right-hand side of \eqref{reladgradbeta}, we have that
\begin{align}
f_i(x^k+td)  \leq &f_i(x_i^k,x_{\neg i}^{k}+td_{\neg i}) +\alpha t \nabla f_i(x_i^k,x_{\neg i}^{k}+td_{\neg i})^Td_i \nonumber \\
&-(1-\alpha)t \theta  \| \nabla f_i(x_i^k,x_{\neg i}^{k}+td_{\neg i})\| \|d_i\| +\frac{t^2L_ k}{2} \|d_i\|^2 \nonumber \\
\leq&f_i(x_i^k,x_{\neg i}^{k}+td_{\neg i}) +\alpha t \nabla f_i(x_i^k,x_{\neg i}^{k}+td_{\neg i})^Td_i \nonumber \\
&-t \left[\frac{(1-\alpha)\theta \lambda_{min}}{(1+\beta)} -\frac{tL_k}{2} \right] \|d_i\|^2 . \nonumber  
\end{align}
Therefore, \eqref{armijo1} and \eqref{armijo2} also hold for $t \leq t'$. Thus, by the way $t$ is updated, we have that $\{t^k\}$ is bounded away from zero.
\halmos 

Let us now consider a class of functions that satisfy the conditions discussed.
\begin{example}\label{exemptafastado}
Consider an equilibrium problem where the functions are given by
\[f_i(x_1,x_2) = \psi_i(x_i) + \nu_i(x_i)^T x_{\neg i}, \quad i=1,2,\]
where $\psi_i\colon\R^{n_i} \to \R$,  $\nu_i\colon\R^{n_i} \to \R^{n_{\neg i}}$  are three times continuously differentiable functions. In this case, $\nabla_{x_i} f_i(x)$ is afine with respect to $x_{\neg i}$. Thus, the Taylor remainder in \eqref{restotaylor1grad} is zero, so \eqref{eq:rmenorbetagrad} holds trivially. Therefore, Proposition \ref{teo:restobetatafastado0} guarantees that $\{t^k\}$ remains bounded away from zero whenever the sequence $\{x^k\}$ is bounded.
\end{example}

Next we present an illustrative example of a non-convex case with a non-zero remainder in which we can also guarantee \eqref{eq:rmenorbetagrad}.
\begin{example}\label{exemptafastado2}
Given $a,b,c>0$, consider the one-dimensional equilibrium problem defined by
$f_i(x) =  2ax_i^2- 2abe^{-cx_i^2} + x_i^4x_{\neg i}^2 $, $i=1,2$. Since $f_i(x) > f_i(x^k_i, x_{\neg i})$ whenever $|x_i| > |x^k_i|$, we have that $\{x^k\} \subset B(0, \|x^0\|)$. Thus, by Propositions \ref{lemah3} and \ref{constantesuniformes}, we have that $\|d\| \leq C_d$, for $t \leq \frac{\lambda_{min}^2}{8\lambda_{max}C_H}$. Moreover, we have that
\[
|\nabla_{x_i} f_i(x)| = \left| 4a x_i \left(1 + bce^{-c x_i^2} + \frac{x_i^2 x_{\neg i}^2}{a} \right) \right| \geq 4a |x_i|
\]
and
\[
|r_i^k(t, d_{\neg i})| = |4 t^2 d_{\neg i}^2 (x_i^k)^3| \leq 4 C_d^2 \|x^0\|^2 |x_i^k|.
\]
Therefore, for any $0\leq\beta <\lambda_{min}/(\lambda_{max}+\lambda_{min})$, if $\|x^0\| \leq \frac{\sqrt{a \beta}}{C_d}$, we have that \eqref{eq:rmenorbetagrad} holds, and consequently, $\{t^k\}$ is bounded away from zero.
\end{example}


Looking from the perspective of traditional optimization problems, one would expect the limitation of $t$ for more general cases, since this occurs on a neighborhood of a point where the gradient is Lipschitz continuous, see \eqref{tafstado0otim}. However, NEPs can behave surprisingly differently than optimization problems. For instance, in \cite{rojas} it is shown that, in constrained problems, it is not always possible to force the KKT residual to be as close to zero as desired on a neighborhood of a solution, unlike the case of optimization problems. In the case of the algorithm proposed in this work, we highlight that the approximation of the gradient done in \eqref{aproxgrad} can demand step-sizes arbitrarily close to zero in order to achieve good accuracy and thus guarantee  descent directions around a point, as is shown in the next example.

\begin{example}\label{novo2}
Consider the equilibrium problem where $f_1(x_1,x_2) := \frac{x_1^2}{2}-2x_1x_2^2$ and
$f_2(x_1,x_2):= \frac{x_2^2}{2}+x_2$. The derivatives are given by 
\begin{itemize}
	\item $\nabla_{x_1} f_1(x)=x_1 - 2x_2^2$, $\nabla^2_{x_1x_1} f_1(x)=1$, 
	$\nabla^2_{x_1x_2} f_1(x)= -4x_2$;
	\item $\nabla_{x_2} f_2(x)=x_2+1$, $\nabla^2_{x_2x_2} f_2(x)=1$ and $\nabla^2_{x_2x_1} f_2(x)=0$. 
\end{itemize}

For an arbitrary $\epsilon>0$, starting from $(x_1^k,x_2^k):=(-\epsilon,\sqrt{\epsilon})$ and using $H_1^k= \nabla^2_{x_1x_1}f_1(x^k)=1$ and $H_2^k=\nabla^2_{x_2x_2} f_2(x^k)=1$, we obtain the Newtonian system 
\eqref{NewtonsistKKT2}: 
\begin{equation*} \left(
\begin{array}{cc}
1 & -4t\sqrt{\epsilon} \\[3pt]
0&1
\end{array} \right)
\left(\begin{array}{c}
d_1^{}\\
d_2^{}
\end{array} \right)=
\left(\begin{array}{c}
3\epsilon\\[3pt]
-\sqrt{\epsilon}-1
\end{array} \right),
\end{equation*}
giving $d_1=3\epsilon-4t\sqrt{\epsilon}(1+\sqrt{\epsilon})$ and $d_2=-(1+\sqrt{\epsilon})$. Therefore, we have that 
$\nabla_{x_1}f_1(x^k_1,x_2^k+td_2)
=-\epsilon-2(\sqrt{\epsilon}-t(1+\sqrt{\epsilon}))^2<0.$
Hence, in order for the direction $d_1$ to be a descent direction, we need $d_1>0$, that is,
$t < \frac{3 \sqrt{\epsilon}}{4(1+\sqrt{\epsilon})}$. Therefore, we see that the step-size may need to be arbitrarily reduced in order for the step to be accepted. This example shows that the situation is not analogous to the case of optimization problems and the search for conditions guaranteeing that $t^k$ is bounded away from zero may not be as simple as in the case of optimization problems.

\end{example}


Also, unlike in optimization analyses, when the sequence $\{x^k\}$ is not necessarily convergent, one may not expect to prove stationarity of its limit points.  However, in the next proposition we show that under certain conditions, if there is some indication of asymptotic stationarity for the first player, then one can infer that either the other player is also stationary or its function is unbounded from below. That is, even in the absense of convergence of the sequence, and under reasonable assumptions, the algorithm finds stationary points when the problem is not ill-posed. 
Clearly, the role of each player may be exchanged in the statement of this result.

\begin{proposition}\label{sumariovelho}
	Consider $\{x^k\}$ and $\{t^k\}$ given by Algorithm \ref{Algorithm 1} and suppose that  there exists a positive real number $\overline{t}$ such that $t^k > \overline{t}$ for all $k \in \mathbb{N}$. 
 \begin{itemize}
     \item[i)] If there exists $k_0$ such that $d_2^k=0$ for all $k \geq k_0$ then either $\{f_1(x^k)\}$ is unbounded from below, or $\lim_{k\to+\infty}g_1^k =0$.
     \item[ii)] If Assumption \ref{Asssump:RestoGrad} holds with $C_R^k$ and the mixed gradients and Hessians are  uniformly bounded and $\sum_{k=1}^{\infty} \|d_2^k\|<+\infty$, then either $\{f_1(x^k)\}$ is unbounded from below, or there exists an infinite subset $\N' \subset \N$ such that $\lim_{k \in \N'}g_1^k =0$.
 \end{itemize}  
\end{proposition}

{\bf Proof:} By conditions \eqref{armijo1} and \eqref{armijot1}, we have that
\begin{align*}	f_1(x_1^{k+1},x_2^{k+1})&\leq f_1(x_1^{k},x_2^{k+1}) +t^k\alpha\nabla_{x_1}f_1(x_1^{k},x_2^{k+1})^Td_1^{k} \nonumber \\
	& \leq f_1(x_1^{k},x_2^{k+1})-t^k\alpha\theta\|\nabla_{x_1}f_1(x_1^{k},x_2^{k+1})\|\cdot\|d_1^{k}\|. 
\end{align*}

Summing and subtracting $f(x^k)$, using that $\bar{t} \leq t^k \leq 1$ and the mean value theorem, we have for some $\xi_k \in [0,1]$ that
\begin{align}
	f_1(x^{k+1}) - f_1(x^{k}) &\leq f_1(x_1^{k},x_2^{k+1}) - f_1(x_1^{k},x_2^{k}) -t^k\alpha\theta\|\nabla_{x_1}f_1(x_1^{k},x_2^{k+1})\|\cdot\|d_1^{k}\| \nonumber \\
	&\leq t^k\nabla_{x_2}f_1(x_1^k,x_2^k +\xi_kt^kd_2^{k})^Td_2^{k}-t^k\alpha\theta\|\nabla_{x_1}f_1(x_1^{k},x_2^{k+1})\|\cdot\|d_1^{k}\| \nonumber \\
	&\leq  \|\nabla_{x_2}f_1(x_1^k,x_2^k +\xi_kt^kd_2^{k})\| \|d_2^{k}\|-\bar{t}\alpha\theta\|\nabla_{x_1}f_1(x_1^{k},x_2^{k+1})\|\cdot\|d_1^{k}\|. \label{mediaria}
	\end{align}
Using \eqref{armijod1} and a telescoping sum on \eqref{mediaria}, if $d_2^k=0$ for $k>k_0$ then
\begin{equation*}
f_1(x^{k+1})\leq f_1(x^{k_0})-\bar{t}\alpha\theta\gamma \sum_{j=k_0}^n \|g_1^k\|^2.
\end{equation*}
So, if $\{f_1(x^{k})\}$ is bounded, we must have that $\lim_{k \to \infty} g_1^k=0$. 

 Let us consider now the second case. Since the mixed gradients, Hessians and $C_R^k$ are bounded by an universal constant, and  since $d_2^k \rightarrow 0$, Assumption \ref{Asssump:RestoGrad} ensures that given $\epsilon>0$ there exist $C_g>0$ and $k_1\in \mathbb{N}$ such that if $\|g_1^k\| > \epsilon$ and $k>k_1$ then 
\begin{align*}
& \|\nabla_{x_2}f_1(x_1^k,x_2^k +\xi_kt^kd_2^{k})\| \leq 2 C_g,\\
& \|\nabla_{x_1}f_1(x_1^k,x_2^k +t^kd_2^{k})\| \leq 2\|g_1^k\|,\\
& \|d_2^k\| \leq \frac{\overline{t}\alpha\theta \gamma \epsilon^2}{8 C_g } \leq \frac{\overline{t}\alpha\theta \gamma }{8 C_g }\|g_1^k\|^2.
\end{align*}
So, in this situation we have by \eqref{mediaria} and  \eqref{armijod1} that
\begin{align}
    f_1(x^{k+1}) - f_1(x^{k}) &\leq 2 C_g\cdot\|d_2^{k}\| -\frac{\overline{t}\alpha\theta}{2} \|g_1^k\|\cdot\|d_1^{k}\| \nonumber\\
    & \leq \frac{\overline{t}\alpha\theta \gamma}{4 }\|g_1^k\|^2-\frac{\overline{t}\alpha\theta \gamma}{2 }\|g_1^k\|^2\leq -\frac{\overline{t}\alpha\theta \gamma \epsilon^2}{4 } \label{mediariavelha}.
\end{align}

If \eqref{mediariavelha} holds for all $k$ large enough we have that $\{f_1(x^{k})\}$ is unbounded. Thus, assuming that it is bounded, there exists an infinite subset $\N' \subset \N$ such that $\lim_{k \in \N'}g_1^k =0$. Using the first equation of \eqref{NewtonsistKKT2}, the boundedness of the mixed Hessians, Assumption \ref{Assump:Huniforposdeflim}, and the fact that $d_2^k \to 0$, we conclude that $\lim_{k \in \N'}d_1^k =0$ and so the result follows from Theorem \ref{pro:conv}.
\halmos

 As a consequence of Proposition \ref{sumariovelho},  item i) states that, considering that $f_1$ is bounded, every limit point of $\{x^k\}$ is stationary for the NEP \eqref{PJ1}, while item~ii) states that, if the sequence lies in a compact set, there exists a limit point that is stationary. The results are, in some sense, a generalization of what happens when Algorithm \ref{Algorithm 1} is applied to a problem where optimality for one of the player is always satisfied, and only optimality of the other player is to be sought. In this case, $t^k$ is always bounded away from zero and the Newton direction has a norm greater than a fraction of the gradient. Therefore, our result recovers the standard  result from Newton's method that either the function is unbounded from below or all  limit points are stationary. 
 
Next we show that, when the objective functions are strictly convex quadratic functions, the full step with $t=1$ is always accepted by the line-search and one iteration is enough to solve the problem. This is similar to what is known for the standard Newton method for optimization. %

\begin{proposition}\label{quadratico}
Let $A_1 \in \R^{n_1 \times n_1 }$, $A_2 \in \R^{n_2 \times n_2}$, $B_1 \in \R^{n_1 \times n_2}$,  $B_2 \in \R^{n_2 \times n_1}$, $c_1 \in \R^{n_1}$, and $c_2 \in \R^{n_2}$ with $A_1$ and $A_2$ being positive definite matrices, with eigenvalues in $[\lambda_{min}, \lambda_{max}]$, $\lambda_{min}>0$. Consider the NEP defined by the functions $f_1(x_1,x_2)=\frac{1}{2} (x_1)^TA_1 x_1+(B_1x_2-c_1)^Tx_1$ and $f_2(x_1,x_2)=\frac{1}{2} (x_2)^TA_2 x_2+(B_2x_1-c_2)^Tx_2$. If $\alpha \leq \frac{1}{2}$, $\gamma \leq \frac{1}{\lambda_{max}}$, 
 $\theta \leq \frac{\lambda_{min}}{\lambda_{max}}$, $H_1^0=A_1$, $H_2^0=A_2$, and the matrix $\left(\begin{array}{cc}A_1 & B_1\\B_2 &A_2\end{array} \right)$ is non-singular, then Algorithm \ref{Algorithm 1} finds a solution in a single iteration without performing backtracking; ie, the solution of system \eqref{NewtonsistKKT2} with $t=1$ satisfies all inequalities (\ref{armijo1})-(\ref{armijod2}) 
 and the new iterate $x^1$ is the solution of the problem.

\end{proposition}

{\bf Proof:} We prove the result for $i=1$. For $i=2$ the result is similar. Since $f_1(x)$ is quadratic with respect to $x_1$ and linear with respect to $x_{2}$, we have that $\nabla_{x_1}f_1(x)$ is linear. Thus, given an arbitrary $x$ and $d$,
\begin{align}	 \nabla_{x_1}f_1(x_1+d_1,x_{2})  = \nabla_{x_1}f_1(x) + A_1d_1,\quad \text{ and }\label{gradlin1}\\
	 \nabla_{x_1}f_1(x_1,x_{2}+d_{2})  = \nabla_{x_1}f_1(x) + B_1d_{2}.\label{gradlin2}
\end{align}
So, using \eqref{gradlin2} in \eqref{NewtonsistKKT2} we obtain that
\begin{equation}\label{NetonQk+1}
	A_1 d_1=-\nabla_{x_1}f_1(x_1^0,x_{2}^0+d_{2}). 
\end{equation}
Then, by \eqref{gradlin1}, $\nabla_{x_1}f_1(x^0+d)=0$. Since $f_i(\cdot,x_{\neg i})$, $i \in \{1,2\}$, are convex, we have that $x^0+d$ is a solution of the NEP.

Next, we show that the inequalities (\ref{armijo1})-(\ref{armijod2}) hold, which means that the point $x^1:=x^0+d$ is accepted by the algorithm. By \eqref{NetonQk+1} we have that
$$\|\nabla_{x_1}f_1(x_1^0,x_{2}^0+d_{2})\| \leq \|A_1\| \|d_1\| \leq \lambda_{max} \|d_1\|.$$
So we obtain (\ref{armijod1}) and (\ref{armijod2}) for $\gamma \leq \frac{1}{\lambda_{max}}$.  As a consequence of this and \eqref{NetonQk+1}
$$-d_1^T \nabla_{x_1}f_1(x_1^0,x_{2}^0+d_{2})  =d_1^TA_1d_1 \geq \lambda_{min} \|d_1\|^2 
 \geq  \lambda_{min} \gamma \|d_1\| \|\nabla_{x_1}f_1(x_1^0,x_{2}^0+d_{2})\|.$$
Thus, (\ref{armijot1}) and (\ref{armijot2}) hold. Finally, to prove \eqref{armijo1} and \eqref{armijo2}, we use that $f_1(\cdot,x_{2}^0+d_{2})$ is  quadratic and \eqref{NetonQk+1} to obtain
\begin{align*}
	f_1(x+d)&=  f_1(x_1^0,x_{2}^0+d_{2}) + \nabla_{x_1} f_1(x_1^0,x_{2}^0+d_{2})^Td_1 + \frac{1}{2}d_1^TA_1d_1\\
	&=f_1(x_1^0,x_{2}^0+d_{2}) +  \frac{1}{2}d_1^T\nabla_{x_1} f_1(x_1^0,x_{2}^0+d_{2}).
\end{align*}
\halmos

Note that the conditions concerning the strictly convex quadratic NEP used in Proposition \ref{quadratico} are equivalent to requiring the problem to have a unique solution. This hypothesis is commonly seen or derived from convergence hypotheses in the study of many best response algorithms. An operator frequently employed in the analysis of such algorithms is $\bar{H}:\mathbb{R}^{2n_1+n_2} \to \mathbb{R}^{n_1 \times n_1}$, defined as $
\bar{H} (x_1',x_1'',x_2) :=\nabla^2_{x_1 x_1}f_1 (x_1',x_2)^{-1}\nabla^2_{x_1 x_2}f_1(x_1',x_2)\nabla^ 2_{x_2 x_2}f_2(x_1'',x_2)^{-1}\nabla^2_{x_2 x_1}f_2(x_1'',x_2).$

Using the definition of best response functions given in \eqref{melhorresposta}, we define their composition as $v(x_1):=b_1(b_2(x_1))$. Various connections between $\bar{H}$ and $v$ are well-established, one of which is that the Jacobian of $v(x_1)$ is equal to $\bar{H}(v(x_1),x_1,b_2(x_1))$. For further details, refer to \cite{caruso}, where the authors consider games in arbitrary Hilbert spaces. In \cite{LiBasar}, the authors show that if the spectral radius of $\bar{H}(v(x_1),x_1,b_2(x_1))$ is less than one for all $x_1 \in \mathbb{R}^{n_1}$, then the NEP \eqref{PJ1} has a unique solution, and both the Jacobi and Gauss-Seidel methods converge to the solution from any starting point.

In \cite{Basar}, the author explores various relaxations of the best response methods, combining the current decision with that provided by the traditional iterative scheme. For instance, a relaxation for the Jacobi method could be expressed as $x_i^{k+1}= \eta_ix_i^k+(1-\eta_i)b_i(x_{\neg i}^k),$
where $\eta_i \in \mathbb{R}$ and $i = 1,2$. The study indicates that  using certain possibly non-null values of $\eta_i$ can expand the range of quadratic problems in which the method converges, as well as improve its convergence rate. Similar results are attained for other relaxations of both Jacobi and Gauss-Seidel methods, with certain specific situations in one or two variables explored in greater detail.

In \cite{ref2}, the authors make significant advances in studying relaxations of the type $x_1^{k+1}:= \eta x_1^k+(1-\eta) b_1(b_2(x_1^k)).$
The presented results encompass possibly non-quadratic problems in arbitrary Hilbert spaces, including, as special cases, the results of \cite{LiBasar,Basar,caruso}, among others. To this end, the authors delineate a non-restrictive class of games, termed ratio-bounded games. For finite dimension strategy sets, this class comprises games where there exist $\bar{\alpha},\bar{\beta} \in \mathbb{R}$ such that
$\bar{\alpha} \leq \frac{y^T\bar{H}(x_1',x_1'',x_2)y}{\|y\|^2} \leq \bar{\beta}$
for all $x_1', x_1'',y \in \mathbb{R}^{n_1}$, $y \neq 0$, and $x_2 \in \mathbb{R}^{n_2}$. Accordingly, with $\lambda$ denoting the Lipschitz constant of $v(\cdot)$, the authors define four subclasses of problems based on the values of $\bar{\alpha}$, $\bar{\beta}$, and $\lambda$. For each subclass, they study the existence and uniqueness of the solution, define the range of $\eta$ for which the algorithm converges, and determine its optimal value to achieve the best convergence rate.

It is noteworthy that Proposition \ref{quadratico} not only shows convergence but also ensures it in a single iteration, requiring only the existence and uniqueness of equilibrium in strictly convex quadratic problems. However, we do not believe that our result is directly comparable with those associated with best-response algorithms, as our iteration involves solving a problem of the same class of the original problem (quadratic NEP). Namely, our algorithm directly tackles the $(n_1+n_2) \times (n_1+n_2)$ linear system resulting from the KKT conditions, while in best response methods, one tackles a $n_1 \times n_1$ and a $n_2 \times n_2$ linear system per iteration. In practice, this difference implies that each player needs complete knowledge of the structure of the entire game, which can be quite restrictive. Conversely, in best-response algorithms, each player can base their decision solely on the elements of their own problem, merely observing the other player's decision.

Regarding Yuan's method proposed in \cite{yuan}, while the quadratic model of player $i$ can be constructed solely with $x^k_{\neg i}$, the values of $f_{\neg i}(x^k)$ and $\nabla_{x_{\neg i}} f_{\neg i}(x^k)$ must be shared to determine the trust region radius. Moreover, the positive definiteness of $\left(\begin{array}{cc}A_1 & B_1\\B_2 &A_2\end{array} \right)$ and certain algorithmic assumptions, including the sequence lying in a compact set and an assumption that the step obtained in the minimization of the quadratic model is better than the Cauchy step, are presumed. Nevertheless, the result only ensures the existence of a limit point that is the solution, not guaranteeing the convergence of the entire generated sequence nor establishing results on the speed of the algorithm.

We conclude our convergence results by showing that we can estimate a linear convergence rate when close to a solution, under certain circumstances. For this, for any vector $x \in \mathbb{R}^{n_1+n_2}$ we will denote 
\[H (x) :=\left[\begin{array}{cc}\nabla^2_{x_1x_1}f_1(x) & \nabla^2_{x_1x_2}f_1(x) \\ \nabla^2_{x_2x_1}f_2(x) & \nabla^2_{x_2x_2}f_2(x)\end{array}\right]  \text{ and } 
E (x) :=\left[\begin{array}{cc}0& \nabla^2_{x_1x_2}f_1(x) \\ \nabla^2_{x_2x_1}f_2(x) & 0\end{array}\right].\]

Additionally, we are going to require some control on the Taylor approximation of the gradients around the solution, as established next:

\begin{assump}\label{Asssump:RestoGradSol}
Given a solution $x^* \in \R^{n_1+n_2}$ of the NEP \eqref{PJ1}, this assumption is satisfied at $x^*$ if there exist $\hat{\delta}_*>0,C_R^*>0$, and a function $r\colon\R^{n_1+n_2} \to \R$ such that if $x \in B(x^*,\hat{\delta}_*)$ then \begin{equation}\label{restotaylor1gradSol}
	g(x) = g(x^*)+H(x^*)(x-x^*)+r(x),
\end{equation}
where
	$\| r(x)\| \leq C_R^* \|x-x^*\|^2$.
\end{assump}

\begin{proposition}\label{taxa}
Suppose that $x^*$ is a solution of the NEP \eqref{PJ1} such that Assumption \ref{Asssump:RestoGradSol} holds and $\nabla^2_{x_1x_1}f_1(x^*)$ and $\nabla^2_{x_2x_2}f_2(x^*)$ have all their eigenvalues greater than $\lambda_{min}>0$. Suppose also that $H(x^*)$ is non-singular and that there exist $\delta>0$ and $\hat{t} \in (0,1)$ such that the iteration is well defined with $H_i^k=\nabla^2_{x_ix_i}f_i(x^k)$, for $i =1,2$, and $t^k \geq \hat{t}$ whenever $x^k \in B(x^*,\delta)$, and 
\begin{equation}\label{limitaErroHE}
    \|E(x^*)\| \leq 
    \frac{1}{32\|H(x^*)^{-1}\|}. 
\end{equation} Then there exists a neighborhood of $x^*$ such that if $x^k$ is in this neighborhood, then
\begin{equation}
    \|x^{k+1}-x^*\| \leq \left(1-\frac{\hat{t}}{2}\right) \|x^{k}-x^*\|. \label{teseconvlinear}
\end{equation}

\end{proposition}

{\bf Proof:} Using that $g^*:=g(x^*)=0$, by the definition of $x^{k+1}$, we have that
\begin{equation}
   \|x^{k+1}-x^*\| = \|x^k - x^* - t^k{\cal H}_{t^k,k}^{-1}(x^k)(g^k - g^*)  \|.
\label{primeiranorma}  
\end{equation}
On the other hand, summing and subtracting terms in \eqref{restotaylor1gradSol} for $x=x^k$ we have that
\begin{equation}
    g^k - g^* = {\cal H}_{t^k,k}(x^k-x^*) + \left[H(x^*)- {\cal H}_{t^k,k}\right](x^k-x^*) + r(x^k), \label{ges} 
\end{equation}
with $\|r(x^k)\|\leq C_R^*\|x^k-x^*\|^2$.

Hence, substituting \eqref{ges} in \eqref{primeiranorma} we have
\begin{align}
     \|x^{k+1}-x^*\| &= \left\| (1- t^k)(x^k-x^*)- t^k{\cal H}_{t^k,k}^{-1} \left[(H(x^k)-  {\cal H}_{t^k,k})(x^k-x^*) +r(x^k)\right] \right\| \nonumber\\
    &\leq \left(1-t^k+t^k\|{\cal H}_{t^k,k}^{-1}\|(\|H(x^k)-  {\cal H}_{t,k})\|+C_R^*\|x^k-x^*\|)\right)\|x^k-x^*\|. 
    \label{EqauxConvLocal}
\end{align}
In order to proceed with the limitation of the right-hand side aiming for \eqref{teseconvlinear}, first notice that the continuity of the second derivatives and the fact that $H(x^*)$ is non-singular together with \eqref{limitaErroHE} guarantee that there exists $\hat{\delta}$ such that 
\begin{eqnarray}
    H(x^k) \text{ is invertible and } \|H(x^k)^{-1}\| \leq 2\|H(x^*)^{-1}\|, \label{invH1k}\\
    \|E(x^k)\| \leq \max\left\{2\|E(x^*)\|,\frac{1}{16\|H(x^*)^{-1}\|} \right\}, \label{limE}
\end{eqnarray}
whenever $\|x^k-x^*\| \leq \bar{\delta}:=\min\{\delta,\hat{\delta}\}$.

On the other hand, by the definition of $M_1^k$ and $M_2^k$ and \eqref{limE}, if $x^k \in B(x^*,\bar{\delta})$,
\begin{align}
 \|H(x^k)-{\cal H}_{t^k,k}\| &= \left\| \left[\begin{array}{cc}0& \nabla^2_{x_1x_2}f_1(x^k)-t^kM_1^k \\ \nabla^2_{x_2x_1}f_2(x^k)-t^kM_2^k & 0\end{array}\right]\right\| \nonumber \\
 &\leq \|E(x^k) \| \leq \frac{1}{16\|H(x^*)^{-1}\|}. \label{limitaHmenosHt}
\end{align}
   Then, combining \eqref{invH1k} and \eqref{limitaHmenosHt},  
   \begin{align*}
       \|H(x^k)^{-1}[H(x^k)-{\cal H}_{t^k,k}]\| &\leq \|H(x^k)^{-1}\| \|[H(x^k)-{\cal H}_{t^k,k}]\|\\ & \leq 2\|H(x^*)^{-1}\| \|E(x^k) \| <\frac{1}{2}.
   \end{align*}
Hence, by Banach's Lemma (see Theorem 2.3.4 in \cite{golub}) and \eqref{invH1k}, 
\begin{equation}\label{limHtk}
    \|{\cal H}_{t^k,k}^{-1}\| \leq 2\|H(x^k)^{-1}\| \leq 4 \|H(x^*)^{-1}\|.
\end{equation}

Thus, if $\|x^k-x^*\| \leq \min\left\{\bar{\delta}, \frac{1}{16 C_R^*\|H(x^*)^{-1}\|}\right\}$ we have by \eqref{limHtk} and \eqref{limitaHmenosHt} that
\begin{align}
    \|{\cal H}_{t^k,k}^{-1}\|&(\|H(x^k)-  {\cal H}_{t,k})\|+C_R^*\|x^k-x^*\|) \nonumber\\
    &\leq 4 \|H(x^*)^{-1}\| \left(\frac{1}{16 \|H(x^*)^{-1}\|}+\frac{C_R^*}{16 C_R^*\|H(x^*)^{-1}\|}\right)=\frac{1}{2}. \label{ultimaConvLocal}
\end{align}
Therefore, combining \eqref{EqauxConvLocal} and \eqref{ultimaConvLocal} we obtain
\begin{equation*}
\|x^{k+1}-x^*\| \leq \left(1-t^k+t^k \frac{1}{2}\right)\|x^k-x^*\| \leq \left(1-\frac{\hat{t}}{2}\right)\|x^k-x^*\|,
\end{equation*}
which proves the linear convergence rate stated in \eqref{teseconvlinear}.
\halmos

While the assumptions of Proposition \ref{taxa} are mostly associated with smoothness of the functions defining the problem, condition \eqref{limitaErroHE} is related with the diagonal dominance of $H(x^*)$, noticing that the constant $\frac{1}{32}$ used in \eqref{limitaErroHE} can be improved with some small changes in the proof which we opted to omit. Inspired by the results related with the convergence of Newton's method for optimization problems, it would be natural, however, to expect the quadratic local convergence of Algorithm \ref{Algorithm 1}. The key for proving this in optimization problems is first showing that the step-size $t=1$ is always accepted for convex quadratics when $\alpha \leq \frac{1}{2}$. Since smooth functions behave locally very similarly to its Taylor quadratic approximations, one can obtain that if $\alpha <\frac{1}{2}$ then the unitary step is accepted on a neighborhood of a solution where the second order sufficient conditions hold. However, the situation is very different for the case of NEPs where we have that the approximation of the gradient \eqref{aproxgrad} can make it so that the step-size needs to be reduced in order to attain a good approximation accuracy, regardless of how close we are from a well-behaved solution. The next example illustrates this situation, indicating that it is not straightforward to obtain better than linear local convergence rates. In the example we show that arbitrarily close to the unique solution of a convex NEP, the step-size $t=1$ may be rejected.



\begin{example}\label{exnovo}
Consider the equilibrium problem where $f_1(x_1,x_2) := \frac{x_1^2}{2}-2x_1x_2^2$ and
$f_2(x_1,x_2):= \frac{x_2^2}{2}$. This is a convex NEP with the only solution $x^*=(0,0)$.  Starting from $(x_1^k,x_2^k):=(-\epsilon,\sqrt{\epsilon})$, with $\epsilon>0$ arbitrarily small, we have that the Newton direction $(d_1,d_2)$ obtained from \eqref{NewtonsistKKT2} is $d_1=(3-4t)\epsilon$ and $d_2=-\sqrt{\epsilon}$. In this way, $\nabla_{x_1}f_1(x^k_1,x_2^k+td_2)=-\epsilon-2(1-t)^2\epsilon<0$ and so we need to have $t < \frac{3}{4}$ in order to guarantee that $d_1$ is a descent direction for the first player, while fixing the predicted decision of the second player.

\end{example}

We end this section by remarking that it is possible to generalize Algorithm \ref{Algorithm 1} and the results of this paper to consider a number $N>2$ of players. The extension of the algorithm and the well-definiteness results would be natural. The main result,  item ii) of Theorem \ref{pro:conv}, would state that if the sequence converges then the gradient of at least one of the players converges to zero, while all other convergence results would also hold by simply replacing the assumptions on the behavior of player $\neg i$ by a similar requirement for all players different than player $i$. 

\section{Numerical Experiments}\label{sec4}

In this section we present some numerical experiments depicting the behavior of Algorithm \ref{Algorithm 1}. First we test it on six selected illustrative example problems  and then we test the Algorithm on an application on facility location problems. 

\subsection{Illustrative Examples}\label{subsec1}
For the first experiments, we chose the initial point as $x^0 := (-5,1)^T$, an error tolerance of $10^{-4}$ for the norm of the gradients $\|g^k\|$ as stopping criterion, together with a maximum number of iterations of $k=1000$. In Table \ref{table1} we show the convergence point, the error, and the number of iterations for the six selected equilibrium problems to be defined below. We compared Algorithm~\ref{Algorithm 1} (abbreviated by Alg \ref{Algorithm 1} in the tables) with the Jacobi-type trust-region method of \cite{yuan} (abbreviated Yuan), an exact Jacobi-type method that deals with system \eqref{sistKKT} directly at each step (abbreviated Jacobi), and Newton's method for system \eqref{sistKKT} with unitary step, abbreviated Newton. We also display each problem's actual solution for comparison. 

 Jacobi's method for system \eqref{sistKKT} is done as follows: given $(x_1^k,x_2^k)$ we algebrically compute the next iterate $(x_1^{k+1},x_2^{k+1})$ by solving the equation $\nabla_{x_1}f_1(x_1,x_2^k)=0$ for $x_1$ and $\nabla_{x_2}f_2(x_1^k,x_2)=0$ for $x_2$, and then update $(x_1^{k+1},x_2^{k+1}):=(x_1,x_2)$. As for Newton's method, given an iterate $(x_1^k,x_2^k)$, we compute $(d_1^k,d_2^k)$ by solving system \eqref{std} and then compute $(x_1^{k+1},x_2^{k+1})$ as $(x_1^k,x_2^k)+(d_1^k,d_2^k)$.

Yuan's algorithm,  a  Jacobi-type method for NEPs that also uses Newtonian ideas, depends on several parameters, for details see \cite{yuan}. In all our implementations we chose the parameters $\delta_v:=0.01, \tau_v:=1, t_{v,1}:=1$, $G_i=1$, and $\beta_i:=0.5, i=1, 2$. 

For our method, we chose $\alpha:=10^{-6}, \theta=0.01, \gamma = 10^{-6},$ and $\tau=0.99$. Recall that according to Proposition \ref{welldefinedblue}, $\gamma$ is proportional to $\frac{1}{\lambda_{max}}$ and $\theta$ to the ratio $\frac{\lambda_{min}}{4\lambda_{max}}$. Selecting  parameters $\gamma$ and $\theta$ in this way, we enforce a large bounding interval for $\|H_i^k\|$, which guarantees the validity of Assumption \ref{Assump:Huniforposdeflim}. In order to choose the matrices $H_i^k$, we computed the modified Cholesky decomposition of $\nabla^2_{x_ix_i}f_i(x^k)$ as described in \cite{modcholesky}, which resumes to this matrix whenever it is positive definite, and, when it is not, the procedure returns a positive definite matrix by adding a small diagonal perturbation to the matrix.

 All tests were run in Matlab R2017a in an AMD Ryzen 5 2400G with 3.6Ghz graphics and 8Gb RAM processor. Our illustrative problem set is  composed by the  simple one-dimensional NEPs defined by $f_1$ and $f_2$ given in Table \ref{tabproblemas}.
\begin{table}[h]
\caption{Definition of simple problems}
\label{tabproblemas} 
\begin{tabular}{lll}
Example 5.1:\dummylabel{examp1}{5.1}         & $f_1(x):= x_1^2 + x_1x_2 - 5x_1$, & $f_2(x):= 3x_2^2/2 - x_1x_2 - x_2$. \\
Example 5.2:\dummylabel{examp2}{5.2} & $f_1(x):= x_1^2/4 + x_1x_2 - 5x_1$, & $f_2(x):= x_2^2/6 - x_1x_2 - x_2$.\\
Example 5.3:\dummylabel{examp3}{5.3} &$f_1(x):= x_1^2 + x_1x_2 - 5x_1$,&$f_2(x):= -3x_2^2/2 - x_1x_2 - x_2$.\\
Example 5.4:\dummylabel{examp4}{5.4} &$f_1(x) := -x_1(0.6-x_2)$,&$f_2(x) := x_2(0.7-x_1)$.\\
Example 5.5:\dummylabel{exampnovo6}{5.5} &
$f_1(x) := (x_1+x_2)^2/4 + \sin(x_1)$,&$ f_2(x) := (x_1+x_2)^2/4 +\hspace{-1pt} \sin(x_2)$.\\
Example 5.6:\dummylabel{examp5}{5.6}&$f_1(x):= x_1^3x_2^2/3 + x_1^2/2$,& $f_2(x):= x_1^2x_2^3/3 + x_2^2/2$.
\end{tabular}
\end{table}

The first and second examples are defined by positive definite quadratic functions. In the first one, the underlying linear system given by the first-order necessary optimality conditions (which are sufficient in this case) is strictly diagonal dominant which favors convergence of the Jacobi method, while in the second one the spectral radius of the matrix associated with Jacobi's iteration is grater than one, which may hinder its convergence. The third example has no solution, however there exists a non-equilibrium stationary point at $(3.2,-1.4)$. The forth example has zero second-order derivatives, implying that Jacobi's method is undefined. The fifth example has several solutions and non-equilibrium saddle points. For Jacobi's method we solved the equations using Matlab's \texttt{fsolve} function. Finally, in the sixth example, there is an equilibrium point at $(0,0)$ and a non-equilibrium stationary point at $(-1,-1)$.

The results are given in Table \ref{table1}. For each method we report the last point obtained, the norm of $g^k$, 
and the number of iterations used by each method.
\begin{table}[ht]
	\centering

	\caption{Test results}
	\label{table1}
	\resizebox{\columnwidth}{!}
	{
	\centering
	\begin{tabular}{|c|c|c|c|c|} 
	\hline
	\multicolumn{5}{|l|}{$\qquad\qquad$Example \ref{examp1}: Quadratic with matrix with spectral radius $< 1$ - Solution: $(2,1)$ }\\
	\hline
	Method & Alg \ref{Algorithm 1} & Yuan & Jacobi & Newton \\
	\hline
	Point & $(2,1)$ & $(1.99997,0.99999)$ & $(2,1)$ & $(2,1)$ \\
	\hline
	$\|g^k\|$ & 0 & $9.85431\cdot 10^{-5}$ & $0$ & 0 \\
	\hline
	\# iter & 1 & 1 & 2 & 1 \\
	\hline
	\hline
		\multicolumn{5}{|l|}{$\qquad\qquad$Example \ref{examp2}: Quadratic with matrix with spectral radius $>1$ -  Solution: $(0.57142,4.71428)$}\\
		\hline
	Method & Alg \ref{Algorithm 1} & Yuan & Jacobi & Newton \\
	\hline
	Point & $(0.57142,4.71428)$ & $(0.57145,4.71419)$ & $(\infty,\infty)$ & $(0.57142,4.71428)$  \\
	\hline
	$\|g^k\|$ & $2.22009\cdot 10^{-16}$ & $9.71940\cdot 10^{-5}$ & $\infty$ & $2.22009\cdot 10^{-16}$ \\
	\hline
	\# iter & 1 & 158 & 398 & 1 \\
	\hline
	\hline
	\multicolumn{5}{|l|}{$\qquad\qquad$Example \ref{examp3}: Quadratic problem without a solution.}\\
	\hline
	Method & Alg \ref{Algorithm 1} & Yuan & Jacobi & Newton \\
	\hline
	Point & $(-\infty,\infty)$ & $10^{3} \cdot (-0.49401,0.99900)$ & $(3.19999,-1.39999)$ & $(3.2,-1.4)$  \\
	\hline
	$\|g^k\|$ & $\infty$ & $2.50109\cdot 10^{3}$ & $4.28624 \cdot 10^{-5}$ & $1.77615\cdot 10^{-16}$ \\
	\hline
	\# iter & 372 & 1000 & 8 & 1 \\
	\hline
	\hline
	\multicolumn{5}{|l|}{$\qquad\qquad$Example \ref{examp4}: Problem with null Hessians - Solution: $(0.7,0.6)$ }\\	
			\hline
	Method & Alg \ref{Algorithm 1} & Yuan & Jacobi & Newton \\
	\hline
	Point & $(0.7000,0.60001)$ & $10^9 \cdot (1.86823,-0.62812)$ & - & $(0.7,0.6)$  \\
	\hline
	$\|g^k\|$ & $7.00011\cdot 10^{-5}$ & $2.59973\cdot 10^9$ & - & $2.22009\cdot 10^{-16}$ \\
	\hline
	\# iter & 1 & 1000 & - & 1 \\
	\hline
 \hline
	\multicolumn{5}{|l|}{$\qquad\qquad$Example \ref{exampnovo6}: Non-quadratic problem with several equilibria }\\	
			\hline
	Method & Alg \ref{Algorithm 1} & Yuan & Jacobi & Newton \\
	\hline
	Point & $(-2.4026,3.8812)$ & $(-3.14153,3.14153)$ & $(0.73908,0.73908)  $ & $10^3\cdot (-3.07090,3.07719)$  \\
	\hline
	$\|g^k\|$ & $5.87943\cdot 10^{-5}$ & $9.92824\cdot 10^{-5}$ & $2.59534$ & $6.28318$ \\
	\hline
	\# iter & 5 & 19 & 1000 & 1000 \\
	\hline
	\hline
	\multicolumn{5}{|l|}{$\qquad\qquad$Example \ref{examp5}: Cubic problem - Solution: $(0,0)$}\\	
	\hline
	Method & Alg \ref{Algorithm 1} & Yuan & Jacobi & Newton \\
	\hline
	Point &  $10^{-8} \cdot (0.27025,0.27025)$ & $10^{-4} \cdot (0.64551,0.31008)$  & $ (-10^{-154} \cdot 0.74020,-\infty)$ & $ (-1,-1)$  \\
	\hline
	$\|g^k\|$ & $5.40980\cdot 10^{-9}$ & $8.50096\cdot 10^{-5}$ & $\infty$ & $5.04872\cdot 10^{-16}$ \\
	\hline
	\# iter & 9 & 25 & 6 & 7 \\
	\hline

\end{tabular}
	}
	\end{table}
In Examples \ref{examp1} and \ref{examp2} we obtained the expected results: Jacobi managed to converge for the first one but not for the second one, while all the other three methods found the solution. Algorithm \ref{Algorithm 1} and Newton's method converged in one iteration since the problems are quadratic with positive definite Hessians.

In Example \ref{examp3}, which does not have a solution, the exact Jacobi and Newton's methods converged to the stationary point $(3.2,-1.4)$, meaning that they did not aim at solving the actual problem. Our algorithm and Yuan's algorithm, on the other hand, did not converge, as should be expected. Yuan's method performed the maximum number of iterations, while Algorithm \ref{Algorithm 1} found $(-\infty,+\infty)$ in $372$ iterations.

In Example \ref{examp4} the exact Jacobi method is undefined, since the objective functions' Hessians $\nabla_{x_ix_i}^2f_i(x^k)$ are null, meaning that the component-wise systems in \eqref{sistKKT} are undefined. Since the problem is quadratic, Newton's method managed to solve it in one step, as expected. Our method does not use the Newton step here, as the Hessians are not positive definite, thus implying it had to compute a positive definite approximation. Even so, the method still managed to find a solution with the required precision in just one step. We highlight that Yuan's method could not find the solution in this example, due to the fact that Yuan's method relies on the positive definiteness of the Jacobian of the KKT system, which does not occur here. On the other hand, our method managed to reach the solution.

In Example \ref{exampnovo6} there are several equilibrium points due to the nature of the sine function. Both  Algorithm \ref{Algorithm 1} and Yuan's method converge to stationary points that are actually equilibrium. Starting from the aforementioned starting point, neither Jacobi nor Newton's method managed to converge. We point out, however, that in this example Newton's method is very sensible to the starting point. Starting with $(-3,1)$ for instance, yields the point $(-1.7509,1.5709)$ in 5 iterations, which is stationary, but it is a saddle point.

As for Example \ref{examp5}, since the gradients are given by $x_i(x_ix_{\neg i}^2+1)$ for each component $i$, the stationary points can be achieved either if $x_i=0$ or if $x_ix_{\neg i}^2+1=0$. In order to implement the exact Jacobi method for system \eqref{sistKKT}, we need to choose how we are going to satisfy the system directly, either by selecting $x_i=0$ (the equilibrium point) or selecting $x_i$ satisfying $x_ix_{\neg i}^2+1=0$. If we choose the former, Jacobi's method unsurprisingly takes one step. If instead we choose the latter, Jacobi's method behaves as follows: one of the components diverges decreasing the function value indefinitely and the other component manages to find a point that zeroes the respective gradient. On the other hand, Algorithm \ref{Algorithm 1}, Yuan's and Newton's method successfully managed to find a stationary point in this example. However, Newton's method converged to $(-1,-1)$, the closest stationary point,  instead of to the actual solution. Algorithm \ref{Algorithm 1} and Yuan's method, on the other hand, converged to the actual solution.

We highlight that Algorithm \ref{Algorithm 1} never reduced the step-size in order to achieve conditions (\ref{armijo1})-(\ref{armijod2}) in the quadratic Examples \ref{examp1}-\ref{examp4}. This is expected and consistent with Proposition \ref{quadratico} for the strictly convex case. In Examples \ref{exampnovo6} and \ref{examp5} the step-size was reduced at most once for the direction computed to be accepted by the descent criteria (\ref{armijo1})-(\ref{armijod2}). 

We have also run Examples \ref{novo2} and \ref{exnovo} from the mentioned starting points with $\epsilon=10^{-4}$. Algorithm \ref{Algorithm 1} converges to the solution $(0,0)$ in $4$ iterations for Example \ref{exnovo}, where the step-size was only reduced to $0.5$ in the first iteration. As for Example \ref{exnovo}, the step-size is greatly reduced in the first iteration as expected, being only accepted with $t=2^{-8}$. After the third iteration, however, the step-size is always accepted with $t=1$ until convergence to the solution $(2,-1)$ in $5$ iterations.




	\subsection{Facility location problem}
	In this experiment we consider a facility location problem, as first studied in \cite{hotel}. In this problem one must choose the best place to open facilities such as to best serve the costumers.  In general, servers seek to be close to their customers due to the various commercial and logistical advantages that this can entail. Here we consider the problem of two players that want to open a facility each in order to obtain the highest expected return on trading with $N$ customers. A player's chances of serving a certain customer depends on the proximity of their facility with the position of the customer in relation to the competitor's proximity. Denoting by $x_i$ the position of facility $i$ for $i=1,2$, and by $z_j$ the position of client $j$ for $j=1,\dots,N$, let us consider here that the probability of player $i$ serving  client $j$ is $1-\frac{\| x_i - z_j\|^2}{\|x_i - z_j\|^2+\|x_{\neg i} - z_j\|^2}$. Thus, if $b_j^i$ is the profit of player $i$ in serving customer $j$, the function to be minimized by each player to obtain the best expected revenue is
 \begin{equation}\label{fis4}
 f_i(x_i,x_{\neg i}) = \displaystyle \sum_{j=1}^N \frac{b^i_j\|x_i - z_j\|^2}{\|x_i - z_j\|^2+\|x_{\neg i} - z_j\|^2}, \quad i=1, 2,
\end{equation}
resulting in a NEP.

\begin{figure}[ht!]
\begin{center}
    \includegraphics[width=0.350\textwidth]{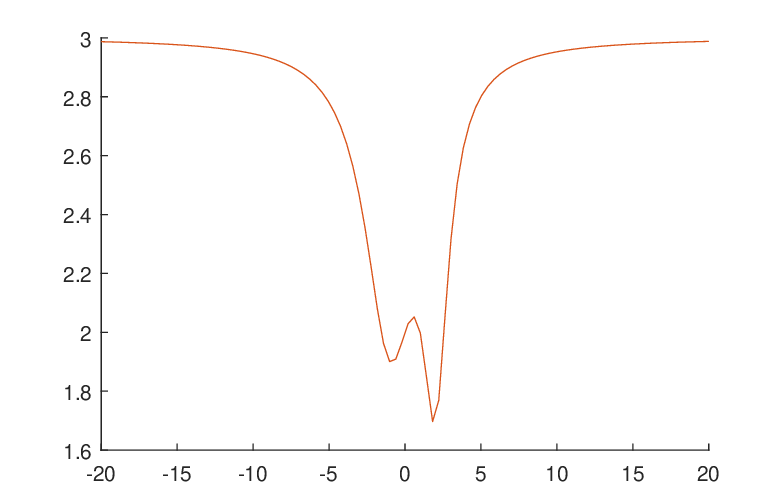}
\end{center}
\caption{Unidimensional objective function for the  facility location problem.}
\label{cortefigsox}
\end{figure}

Notice that the functions are not defined when both facilities are located exactly at the position of one of the customers, a very pathological situation which was never encountered by any of the algorithms we run. 
Anyway, these problems are complex, since the functions are non-convex. In Figure \ref{cortefigsox} we illustrate the objective function of the first player in a one-dimensional problem, with the decision of the second player fixed at $0.915$ and three clients located at $z_1=1, z_2=-1$ and $z_3=3$, with a uniform profit of $1$. For this problem an equilibrium solution is $(1.901,0.915)$, which is found by Algorithm \ref{Algorithm 1}, with the same parameters used  in Subsection \ref{subsec1}, from the starting point $(2,1)$.

\begin{figure}
\centering
\begin{subfigure}{0.496\textwidth}
    \includegraphics[width=\textwidth]{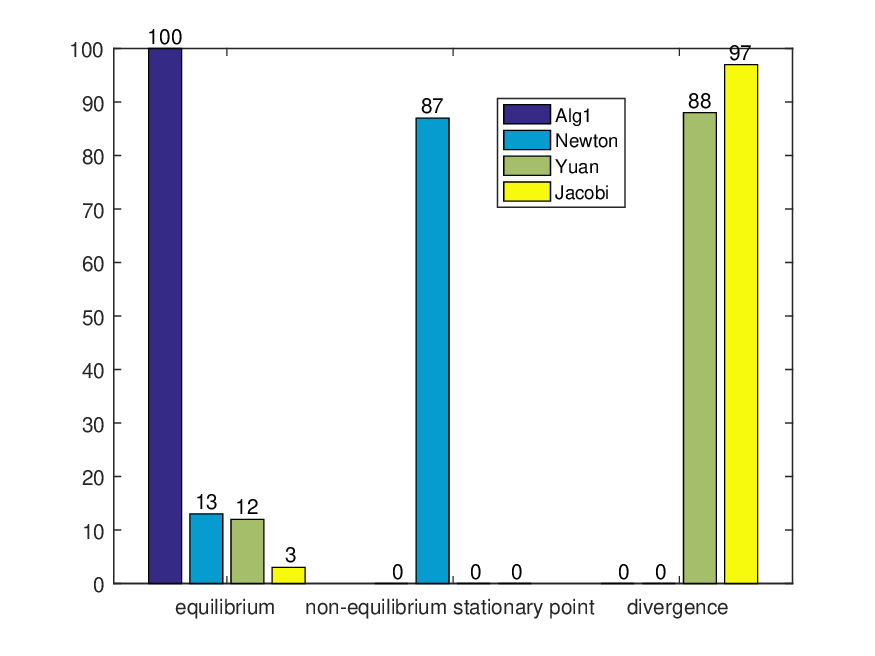}
    \caption{Solution quality for different methods.}
    \label{graficos}
\end{subfigure}
\hfill
\begin{subfigure}{0.496\textwidth}
    \includegraphics[width=\textwidth]{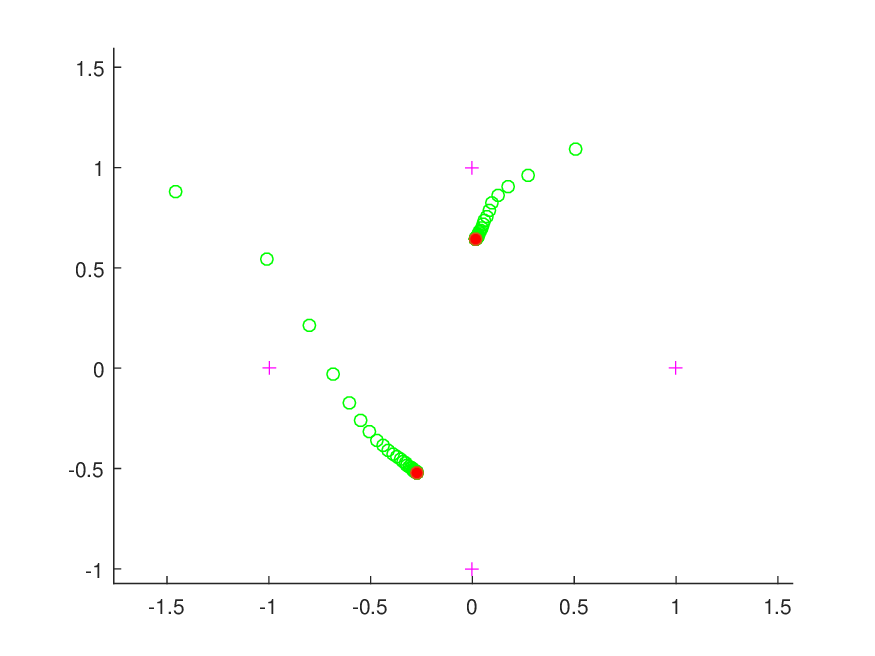}
    \caption{Convergence  trajectory for Alg. \ref{Algorithm 1}.}
    \label{ifdois}
\end{subfigure}
\caption{Numerical  behavior   in the bidimensional facility location problem.}
\label{fig:figures}
\end{figure}

For a more extensive analysis, let us now consider a two-dimensional problem. For our test we selected a problem in the form \eqref{fis4} with the clients' coordinates $z_1=(1,0), z_2=(0,1), z_3=(-1,0)$, and $z_4=(0,-1)$. We chose the profit vector for each player as $b^1=(1,2,1,1)$ and $b^2=(1,2,2,3)$.  We run Algorithm \ref{Algorithm 1} for $100$ random initial points with coordinates generated in the interval $[-2,2]$, comparing with Yuan's, Jacobi and Newton's method. Here, for the Jacobi method we solved problems \eqref{fjac1} 
with $\phi_i = f_i$ using Matlab's \textit{fminunc} solver instead of dealing with system \eqref{sistKKT}. In Figure \ref{graficos} we report the final outcome of each algorithm according to the stopping criterion $\|g^k\|\leq 10^{-6}$. We report simply whether the method converged to an equilibrium point, a non-equilibrium stationary point or if it diverged. We see that Jacobi and Yuan's method perform poorly due to the properties of the Jacobian of the correspondent first-order conditions. Yuan diverged in $88$ out of $100$ runs while Jacobi's method diverged in $97$ runs, reaching the maximum number of iterations. Algorithm \ref{Algorithm 1} converged to a true equilibrium point in all cases, taking on average $36$ iterations and $1.08$ seconds of CPU time, while Newton's method mostly converged to a non-equilibrium stationary point, having attained the true equilibrium point in only $13$ of the $100$ runs. Newton's method, on average, took $19$ iterations and $0.76$ seconds to converge, however, for the $13$ initial points where both algorithms converged to the equilibrium, they were able to perform similarly as Algorithm \ref{Algorithm 1} took $0.82$ seconds to perform $19$ iterations while Newton's method took $0.75$ seconds to perform $16$ iterations, on average.

In Figure \ref{ifdois} we see the behavior of Algorithm \ref{Algorithm 1} for this problem for the initial point $x^0_1=(2,3)$ and $x^0_2=(-3,2)$. The filled circles in the plane represent the optimal facility locations found for the agents, while an empty circle represents the iterates. We see that Algorithm \ref{Algorithm 1} forces the trajectory to the optimal placement. Running Newton's method with this same data yields a sequence that diverges, forcing the second facility to $(+\infty,-\infty)$, yet still stopping declaring success since the gradients tend to zero.

In conclusion, our numerical experiments attest that our method may avoid converging to a non-equilibrium stationary point, differently from methods that deal only with the nonlinear system of equations obtained by the first order necessary optimality conditions as the standard Newton method. Also, differently from Jacobi and Yuan's method, our method is based on a predicted behavior of the other player, which seems to considerably speed up convergence. Moreover, Jacobi and Yuan's method seem to struggle when the Jacobian of the system of equations is not well-behaved, which does not occur with our method.

\section{Final remarks}\label{sec5}


The contribution of this paper is to propose an algorithm for NEPs  with two main characteristics: Firstly, our algorithm addresses the NEP in its entirety, considering its minimization structure and avoiding convergence to some undesired stationary point. Secondly, the iterates of our algorithm are computed in such a way that they can be interpreted as mimicking the behavior of the agents in practice, namely, the iterates computed by the method approximate in some sense the true decision made by the players in different time periods. That is, a new iterate is computed by player one by minimizing a convex quadratic model of their objective function parameterized by a prediction of the action that will be taken by the second player. The prediction is computed simultaneously, assuming that the second player is adopting the same strategy. This gives rise to a Jacobi-type Newton method, globalized by an Armijo-type linesearch, which guarantees that at the new iterate both players will decrease their predicted objectives; this avoids, as much as possible, convergence to undesired stationary points. On the other hand, the method can also be seen as a new interpretation of Newton's method, which is known to be intrinsically linked to fast local convergence.

We belive that both characteristics of our algorithm are new and relevant for the literature of NEPs, especially when non-convexities are considered,  favoring convergence to a true equilibrium rather than simply a stationary point. Moreover, our approach opens the path to new and interesting investigations of the behavior of Newtonian iterates built with the goal of mimicking the behavior of true economic agents. Clearly, this is a starting point for this type of investigation and many aspects still need to be elucidated. For example, it would be interesting to interpret the many parameters of the algorithm with the true agents' behavior. We believe that the synergy between the interpretation of the game's dynamics and the sequence generated by an algorithm can inspire both better techniques for solving NEPs and a better understanding of practical situations associated with models of this type.

In this paper, we focused on the  analysis considering convergence of the iterative sequence when the step size sequence $\{t^k\}$ is bounded away from zero. An important open problem for a better understanding of the algorithm would be to prove whether or not it is possible to guarantee the boundedness away from zero of $\{t^k\}$ in more general cases than those we have discussed. Perhaps adaptations in the algorithm could lead to an affirmative answer to this question.

Another point would be to check whether the convergence study may also incorporate important elements of best response algorithms. For example, we believe that our results have no immediate connection with the concepts of (weighted) potential games \cite{Basar} or ratio-bounded games \cite{ref2}. As the decrease obtained in our iterations considers the function defined by the other player's predicted decision, we cannot guarantee monotonicity of the potential function in consecutive iterations. The class of ratio-bounded games brings very interesting elements about the composition of best response functions, but this composition is not directly associated with our method. However, a good point for future research would be to study whether it is possible to extend such concepts to the new dynamics introduced in our work and then use them in the convergence analysis.

Taking into account that our algorithm converges more often to a true equilibrium, we envision some possible extensions of this work. Our intention is to extend the ideas of this paper to consider such parameterized constrained problems, borrowing ideas from Sequential Quadratic Programming, where one considers a linearization of the constraints in a Newtonian framework. In addition, our unconstrained algorithm may be employed in an Augmented Lagrangian framework for the general constrained problem (see \cite{kanzow16} and \cite{rojas}). In this framework, the constraints are penalized and at each iteration of the algorithm an unconstrained NEP must be solved, and finding true solutions of the subproblems rather than mere stationary points should favor the algorithm in finding true solutions of the original problem.

We also highlight that in a general dynamical system, the interest goes far beyond mere convergence to an equilibrium, since several other interesting phenomena may occur. We expect that subsequent research on this topic should focus on understanding and classifying such phenomena, with the analysis carried out to the algorithms built for mimicking the dynamics. More specifically, we suspect that economic systems behave in practice as a dynamical system, where each player considers the current decision of the other players and a model for predicting their future behavior, before optimizing their strategy. We expect that Jacobi-type methods should capture this behavior, predicting accurately the final outcome of a game, even if no equilibrium is reached. For optimization problems, where a solution represents the objective for the decision maker, it is expected to build a sequence converging to a solution to the problem. In an NEP seen as a dynamic system, the dynamics can converge to an equilibrium state, but other interesting phenomena can be observed, such as orbit, cycling or divergence. In fact, there are situations where an equilibrium point is not even beneficial for players. This happens, for example, in the famous Prisoner's Dilemma and in congestion problems (see \cite{dilemaprisioneiro} and \cite{congestionamento}). In this way, an algorithm that captures the dynamics involved in the game, even if it does not obtain convergence, can  sometimes be even more important than the identification of the solution itself.

Another interesting point to be investigated in the future 
is the connection with dynamic learning, reinforcing decisions that have better objective function values. Since the philosophy of using  best response functions together with predicted decisions for the other player is similar to our approach, our ideas may substitute the standard two-step process of independent decision making and probability distribution updates considered in \cite{gamelearning, leslie}. Furthermore, in stochastic versions of the algorithm, the notions of cycling and orbiting may be more naturally understood as generalized convergence, simplifying the analysis.




\end{document}